\documentstyle[12pt]{article}
\begin{document}

\title{Stable triples, equivariant bundles and \\ dimensional reduction}
\author{Steven B. Bradlow\thanks{Supported in part by a NSF-NATO Postdoctoral
Fellowship
and NSF grant DMS 93--03545} \and Oscar Garc\'{\i}a--Prada}


\maketitle

\let\text\mbox

 \newcommand{\Bbb}{{\bf}}
\newcommand{\bold}{{\bf}}
\renewcommand{\Box}{{}}

 \newcommand{\be}{\begin{equation}}
 \newcommand{\ee}{\end{equation}}
 \newcommand{\scez}{\setcounter{equation}{0}}
 \newfont{\bfc}{cmbsy10 scaled 1200}  
 \newfont{\gl}{eufm10 scaled \magstep1}  
 \newcommand{\dA}{{\Bbb A}}
 \newcommand{\dC}{{\Bbb C}}
 \newcommand{\dE}{{\Bbb E}}
 \newcommand{\dN}{{\Bbb N}}
 \newcommand{\dP}{{\Bbb P}}
 \newcommand{\dR}{{\Bbb R}}
 \newcommand{\dT}{{\Bbb T}}
 \newcommand{\dU}{{\Bbb U}}
 \newcommand{\dZ}{{\Bbb Z}}
 \newcommand{\bA}{{\bf A}}
 \newcommand{\bE}{{\bf E}}
 \newcommand{\bI}{{\bf I}}
 \newcommand{\bV}{{\bf V}}
 \newcommand{\bh}{{\bf h}}
 \newcommand{\cA}{{\cal A}}
 \newcommand{\cB}{{\cal B}}
 \newcommand{\cC}{{\Cal C}}
 \newcommand{\cD}{{\cal D}}
 \newcommand{\cE}{{\cal E}}
 \newcommand{\cF}{{\cal F}}
 \newcommand{\cG}{{\cal G}}
 \newcommand{\cH}{{\cal H}}
 \newcommand{\cI}{{\cal I}}
 \newcommand{\cK}{{\cal K}}
 \newcommand{\cL}{{\cal L}}
 \newcommand{\cM}{{\cal M}}
 \newcommand{\cN}{{\cal N}}
 \newcommand{\cO}{{\cal O}}
 \newcommand{\cP}{{\cal P}}
 \newcommand{\cQ}{{\cal Q}}
 \newcommand{\cR}{{\cal R}}
 \newcommand{\cS}{{\cal S}}
 \newcommand{\cT}{{\cal T}}
 \newcommand{\cU}{{\cal U}}
 \newcommand{\cV}{{\cal V}}
 \newcommand{\cW}{{\cal W}}
 \newcommand{\tA}{{\tilde{A}}}
 \newcommand{\th}{{\tilde{h}}}
 \newcommand{\glg}{\hbox{\gl g}} 
 \newcommand{\glh}{\hbox{\gl h}} 
 \newcommand{\glm}{\hbox{\gl m}} 
 \newcommand{\glu}{\hbox{\gl u}} 
 \newcommand{\glD}{{\hbox{\gl D}}}
 \newcommand{\glE}{{\hbox{\gl E}}}
 \newcommand{\glF}{{\hbox{\gl F}}}
 \newcommand{\glM}{{\hbox{\gl M}}}   
 \newcommand{\glV}{{\hbox{\gl V}}}   
 \newcommand{\hglF}{{\hat{\hbox{\gl F}}}}
 \newcommand{\cglF}{{\check{\hbox{\gl F}}}}
 \newcommand{\ra}{\rightarrow}
 \newcommand{\lra}{\longrightarrow}
 \newcommand{\kahler}{K\"{a}hler}
 \newcommand{\im}{\mathop{{\fam0 Im}}\nolimits}
 \newcommand{\Ker}{\mathop{{\fam0 Ker}}\nolimits}
 \newcommand{\aut}{\mathop{{\fam0 Aut}}\nolimits}
 \newcommand{\End}{\mathop{{\fam0 End}}\nolimits}
 \newcommand{\rank}{\mathop{{\fam0 rank}}\nolimits}
 \newcommand{\ext}{\mathop{{\fam0 Ext}}\nolimits}
 \newcommand{\diff}{\mathop{{\fam0 Diff}}\nolimits}
 \newcommand{\tr}{\mathop{{\fam0 Tr}}\nolimits}
 \newcommand{\Hom}{\mathop{{\fam0 Hom}}\nolimits}
 \newcommand{\vol}{\mathop{{\fam0 Vol}}\nolimits}
 \newcommand{\pic}{\mathop{{\fam0 Pic}}\nolimits}
 \newcommand{\ch}{\mathop{{\fam0 ch}}\nolimits}
 \newcommand{\td}{\mathop{{\fam0 td}}\nolimits}
 \newcommand{\YMH}{\mathop{{\fam0 YMH}}\nolimits}
 \newcommand{\Ch}{\mathop{{\fam0 Ch}}\nolimits}
 \newcommand{\lie}{\mathop{{\fam0 Lie }}\nolimits}
 \newcommand{\dbar}{\overline{\partial}}
 \newtheorem{definition}{Definition}[section]
 \newtheorem{lemma}[definition]{Lemma}
 \newtheorem{prop}[definition]{Proposition}
 \newtheorem{thm}[definition]{Theorem}
 \newtheorem{cor}[definition]{Corollary}
 \newcommand{\pf}{{\em Proof}. }
 \newcommand{\ps}{{p^\ast}}
 \newcommand{\qs}{{q^\ast}}
 \newcommand{\su}{{SU(2)}}
 \newcommand{\xp}{{{X\times\dP^1}}}
 \newcommand{\cod}{{\cO(2)}}
 \newcommand{\sig}{\sigma}
 \newcommand{\tauh}{{\hat{\tau}}}
 \newcommand{\vb}{ vector bundle}
 \newcommand{\vbs}{ vector bundles}
 \newcommand{\hvb}{holomorphic vector bundle}
 \newcommand{\hm}{Hermitian metric}
 \newcommand{\crs}{compact Riemann surface}
 \newcommand{\hvbs}{holomorphic vector bundles}
 \newcommand{\HT}{holomorphic triple}
 \newcommand{\ve}{ vortex equation}
 \newcommand{\ves}{ vortex equations}
 \newcommand{\tve}{$\tau$-vortex equation}
 \newcommand{\tves}{$\tau$-vortex equations}
 \newcommand{\cves}{coupled vortex equations}
 \newcommand{\ctves}{coupled $\tau$-vortex equations}
 \newcommand{\ts}{$\tau$-stable}
 \newcommand{\ths}{$\tauh$-stable}
 \newcommand{\is}{invariantly stable}
 \newcommand{\sui}{$\su$-invariant}
 \newcommand{\sue}{$\su$-equivariant}
 \newcommand{\suis}{$\su$-invariantly stable}
 \newcommand{\hym}{Hermitian--Yang--Mills}
 \newcommand{\ymh}{Yang--Mills--Higgs}
 \newcommand{\he}{Hermitian--Einstein}
 \newcommand{\hymh}{Hermitian--Yang--Mills--Higgs}
 \newcommand{\CMP}{Comm. Math. Phys.}
 \newcommand{\JDG}{J. Diff. Geom.}
 \newcommand{\qed}{\hfill$\Box$}
 \newcommand{\duy}{Donaldson, Uhlenbeck and Yau }
 \newcommand{\kah}{\kahler}
 \newcommand{\cms}{{\cm_\sig}}
 \newcommand{\mt}{{\glM_\tau}}
 \newcommand{\Mt}{{M_\tau}}
 \newcommand{\os}{{\omega_\sig}}
 \newcommand{\TT}{{\theta_\tau}}
 \newcommand{\ms}{{\mu_\sig}}
 \newcommand{\ma}{{\mu_\alpha}}
 \newcommand{\tri}{{(E_1,E_2,\Phi)}}
 \newcommand{\dtri}{{(E_2^\ast,E_1^\ast,\Phi^\ast)}}
 \newcommand{\stri}{{(E_1',E_2',\Phi')}}
 \newcommand{\extn}{{0\lra \ps E_1\lra F\lra \ps E_2\otimes\qs\cod\lra 0}}
 \newcommand{\sextn}{{0\lra \ps E_1'\lra F'\lra \ps E_2'\otimes\qs\cod\lra 0}}
 \newcommand{\XxP}{X\times\Bbb P^1}
 \newcommand{\mm}{missing material to be added}
 \newcommand{\ci}{C^\infty}

 \def \aspace{\Cal C \times \kform 0E}
 \def \bdl{E\longrightarrow X}
 \def \BoB{\Cal B_0/\Cal B}
 \def \Bt{\Cal B_{\tau}}
 \def \curvdh{F_{\dbare,H}}
 \def \Cx{C^{\dbare}_{\phi}}
 \def \dbar{\overline{\partial}}
 \def \dbare{\overline{\partial}_E}
 \def \dbf{\dbare,\Phi}
 \def \G{\frak G}
 \def \Gc{\frak G_\Bbb C}
 \def \Gco{\frak G_{\Bbb C 0}}
 \def \Gcp{\frak G_{\Bbb C p}}
 \def \Hm{Hermitian\ }
 \def \HE{Hermitian-Einstein\ }
 \def \HP{holomorphic pair\ }
 \def \kform#1#2{\varOmega^{#1}(#2)}
 \def \klerform{\frac{1}{4\pi^2}\hat\Omega}
 \def \lieG{\frak g}
 \def \lieGc{\lieG_{\Bbb C}}
 \def \mmapt{\varPsi_{\tau}}
 \def \muM{\mu_M}
 \def \mum{\mu_m (\phi)}
 \def \pqform#1#2#3{\varOmega^{#1,#2}(#3)}
 \def \s{\Cal C_{s}}
 \def \ss{\Cal C_{ss}}
 \def \ssp{\Cal H_{ss}}
 \def \St{\Cal S_{\tau}}
 \def \tV{$\tau$-Vortex\ }

 \def \Vo{\Cal V_0}
 \def \Vog{\Vo/\Gc}
 \def \Vt{\Cal V_{\tau}}
 \def \Vtg{\Vt/\Gc}
 \def \zset{\varPsi^{-1}_{\tau}(0)}
 \def \zsetg{\zset/\G}

 \def \dtri{{(E_2^*,E_1^*,\Phi^*)}}
 \def \stri{{(E_1',E_2',\Phi')}}
 \def \tri{{(E_1,E_2,\Phi)}}
 \def \TT{\theta_\tau}
 \def \XxP{X\times\Bbb P^1}

\section{Introduction}\label{intro}

 The Hitchin-Kobayashi correspondence between stable bundles and
 solutions to the \linebreak
Hermitian-Einstein equations allows one to apply
 analytic methods to the study of stable bundles.  One such analytic
 technique, which has not yet been much exploited, is that of
 dimensional reduction.  This is a useful tool for studying certain
 special solutions to partial differential equations; in particular it
 is useful for studying solutions which are invariant under the action
 of some symmetry group. When applied to the Hermitian-Einstein
 equations, it thus provides a way of looking at holomorphic bundle
 structures which are both stable and invariant under some group action on the
 bundle, i.e. of looking at {\em equivariant stable bundles}.

 The main idea in dimensional reduction is the following: Suppose we
 have a partial differential equation defined on a space which has a
 symmetry, i.e. which supports some group action.  Then by integrating
 over the group orbits, any solution which is invariant under the
 group action becomes an object defined in a space of lower dimension than
 that of the original setting for the general solutions.  This lower
 dimensional space is the orbit space of the group action, and in that
 space the special solutions to the original equations  can be
 re-interpreted as ordinary solutions to some new set of equations.

 In particular, suppose that we start with a 4-manifold, $M$, a Lie group $G$\
which acts
 on it, and a complex vector bundle $E$\ to which this action lifts.
 In such a situation, there can be equivariant
 solutions to the Hermitian-Einstein equations, and dimensional
 reduction can be applied. The orbit space $E/G$\ will be a new bundle over
$M/G$,
 the orbit space of the group action on the original 4-manifold.  The
 equivariant solutions to the original Hermitian-Einstein
 equations will be solutions to a new set of equations on the bundle $E/G\lra
M/G$.
For example, the equations introduced by Hitchin in \cite{H}, namely the
Anti-Self-Dual equations on a Riemann surface, can be viewed as arising in this
way.

A special case of the above situation occurs when the 4-manifold $M$\ is a
complex surface, and the orbit space $M/G$\ also admits a complex structure.
Solutions to the Hermitian-Einstein equations then correspond to stable
holomorphic structures on the bundle $E$.  In such a situation, dimensional
reduction acquires an extra, holomorphic interpretation. It results in
information about the equivariant stable bundles on $M$\ being encoded in a
holomorphic interpretation for the dimensionally reduced equations on
$M/G$ .

 These sort of ideas are developed in \cite{GP3}, where they are applied
 to certain $SU(2)$-equivariant bundles over $\xp$.  Here $X$ is a closed
 Riemann surface and the $SU(2)$-action is trivial on $X$ and the
 standard one on $\dP^1$.  In this case, the equivariant holomorphic bundles
over
 $\xp$\ correspond to holomorphic pairs (i.e. bundles plus prescribed
 global sections) over $X$. The dimensional reduction
 of the Hermitian-Einstein equations gives
 the vortex equations, and the stable equivariant bundles on $\xp$\ correspond
(by dimensional reduction) to $\tau$-stable holomorphic pairs on $X$, with
 $\tau$-stability as defined in \cite{B2} and \cite{GP3}.

 However not all the  $SU(2)$-equivariant holomorphic
 bundles over $\xp$\ correspond to holomorphic pairs on
 $X$.  In fact those that do form a rather
 restricted subset of the set of all such equivariant bundles.
 A very natural  relaxation of this restriction leads to a
 class of equivariant bundles on $\xp$\ which still
 corresponds to data on the (lower dimensional) space
 $X$, but not necessarily to holomorphic pairs.  What such bundles
 correspond to is a pair of bundles plus a holomorphic homomorphism
 between them.  We call such data a holomorphic triple.

In this paper we undertake a detailed investigation of holomorphic triples over
the
closed Riemann surface $X$.  In particular, we define, in Section 3, a notion
of stability for such objects. We explore the relationship between the
stability of a triple and the stability
of the corresponding equivariant bundle over $\xp$. An important feature of the
definition is that, like in the case of holomorphic pairs, it involves a real
parameter.  This can be traced back to the fact that the definition of
stability for a bundle over $\xp$ depends on  the polarization (choice of
\kahler\ metric) on $\xp$.  We discuss the nature of this parameter, and its
influence on the  properties of the stable triples.  We show for example that

{\em In all cases, with one exception, the parameter in the definition of
triples stability lies in a bounded interval.  The interval is partitioned by a
finite set of non-generic values.}

Our main result is given in Section 4. Loosely speaking, it is that the stable
triples over $X$\ can be considered the dimensional reduction of the stable
equivariant bundles over $\xp$. In other words,

{\em
A holomorphic triple over $X$\  is stable if and only if the corresponding
\linebreak
$SU(2)$-equivariant extension over $\xp$\ is stable}

In \cite{GP3} dimensional reduction is applied to the Hermitian-Einstein
equation on  equivariant bundles over $\xp$.  The result is that on bundles
corresponding to triples over $X$, the equivariant solutions correspond to
solutions to a pair of {\em Coupled Vortex
Equations} on the two bundles in the triple. By combining this result, our
dimensional reduction result for stable bundles, and the Hitchin-Kobayashi
correspondence, we can thus show

{\em There is a Hitchin-Kobayashi correspondence between stability of a triple
and existence of solutions to the Coupled Vortex Equations.}

This is discussed in Section 5.  In Section 6 we discuss the moduli spaces of
stable triples.  By identifying these as fixed point sets of an $SU(2)$-action
on the moduli spaces of stable bundles over $\xp$, we obtain results such as

{\em For fixed value of the stability parameter, the moduli space of stable
triples is a quasi- projective variety. For generic values of the parameter,
and provided the ranks and degrees of the two bundles satisfy a certain
coprimality condition, the moduli space is projective}

In Section 2 we have collected together the basic definitions and background
material that we will need.

 \section{Background and Preliminaries}\label{background}

 \subsection{Basic Definitions}
 Let $X$ be a compact Riemann surface.

\begin{definition}
 A holomorphic triple on $X$ is a triple
 $\tri$ consisting of two holomorphic vector bundles $E_1$ and $E_2$ on $X$
 together with a  homomorphism $\Phi:E_2\lra E_1$, i.e. an element
 $\Phi\in H^0(\Hom(E_2,E_1))$.
\end{definition}

In this paper we will develop a theory of holomorphic
 triples as objects in their own right.  We will also show  how  they arise
from \sue\ holomorphic vector bundles over $\xp$.

 Let $\su$ act on $\xp$ trivially on $X$ and in the standard way on $\dP^1$,
that is we regard $\dP^1$ as the homogeneous space $\su/U(1)$. Let $F$ be a
$\ci$ complex \vb\ over $\xp$.

\begin{definition}
The bundle $F$ is said to be $\su$-{\em equivariant}
 if there is an action of $\su$ on $F$ covering the action on $\xp$.
 Similarly a \hvb\ $F$ is $\su$-{\em equivariant} if it is \sue\ as a $\ci$
bundle
 and in addition the action of $\su$ on $F$ is holomorphic.
\end{definition}

\subsection{Smooth and Holomorphic Equivariant bundles}
 Our main objective in this section  will be the study of \sue\ \hvbs\ on
$\xp$, however, before addressing
 this, we shall  analyse the much easier problem of classifying the \sue\ $\ci$
ones.
  Let $p$ and $q$ be the projections from $\xp$  to the first and second
factors
  respectively.

 \begin{prop}\label{equi-vb}
 Every \sue\ $\ci$ \vb\ $ F$ over $\xp$ can be
 equivariantly decomposed,  uniquely  up to isomorphism, as
 $$
  F=\bigoplus_i  \ps E_i\otimes \qs H^{\otimes n_i},
 $$
 \noindent where
 $E_i$ is a $\ci$ \vb\ over $X$, $H$ is the $\ci$ line bundle over $\dP^1$
 with Chern class 1, and $n_i \in \dZ$ are all different.
 \end{prop}
 \pf
 See \cite[Proposition 3.1]{GP3}.

 We shall describe now \sui\ holomorphic structures on a fixed \sue\ $\ci$ \vb\
over
 $\xp$.   We shall restrict ourselves, however, to the case which is relevant
in connection to
 holomorphic triples.
 Let $E_1$ and $E_2$ be $\ci$ \vbs\ on $X$ and let $H$ be as in Proposition
\ref{equi-vb}.
 Consider the \sue\ $\ci$ \vb\
 \be
 F= \ps E_1\oplus \ps E_2\otimes \qs H^{\otimes 2}, \label{equi-f}
 \ee

 Note that the total space of $\ps E_1$ is $E_1\times \dP^1$, and the action of
 $\su$ that we are considering is trivial on $E_1$ and the  standard one on
$\dP^1$;
 similarly for $\ps E_2$. On the other  hand, recall that the
$\su$-equivariant
 line bundle $H^{\otimes 2}$ over $\dP^1\cong \su/U(1)$  corresponds to the one
 dimensional representations of $U(1)$ given by $e^{i2\alpha}$, i.e.
$H^{\otimes 2}=\su\times_{U(1)}\dC$, where
  $(g,v)\sim (g',v')$ if there is an $e^{i\alpha}\in U(1)$ such that
 $g'=e^{-i\alpha}g$ and $v'=e^{i2\alpha}v$.

 The action of $\su$ on $\su\times\dC$, given by
 $$
 \gamma.(g,v)=(\gamma g,v)\;\;\;\;\;\mbox{for}\;\;\;\gamma\in \su\;\;\;
 \mbox{and}\;\;\;(g,v)\in \su\times \dC,
 $$
  descends to an action on $H^{\otimes 2}$.

 In order to avoid the  introduction of  more notation we shall denote a $\ci$
\vb\
 and the same bundle endowed with a holomorphic structure by the same symbol.
 The distinction will be made  explicit unless it is obvious from the context.

\subsection{Dimensional Reduction of Bundles}
 \begin{prop}\label{equi-hvb}
 Every \sue\ \hvb\ $F$ with  underlying \sue\ $\ci$ structure  given by
(\ref{equi-f})
 is in one-to-one correspondence with a holomorphic extension of the form
 \be
 \extn, \label{extension}
 \ee
 where $E_1$ and $E_2$ are the bundles over $X$ defining (\ref{equi-f})
equipped with
 holomorphic structures.

 Moreover, every such extension is defined by an element
$\Phi\in\Hom(E_2,E_1)$, and is thus
  in one-to-one correspondence with the holomorphic triple $\tri$ on $X$.
 \end{prop}
 \pf
 We shall give here only a brief sketch of the proof (see \cite[Proposition
3.9]{GP3}
 for details). Let $E_1$ and $E_2$ be two \hvbs\ over $X$.
 The extensions over $\xp$  of the form (\ref{extension}) are parametrized by
 $$
 H^1(\xp, \ps (E_1\otimes E_2^\ast)\otimes \qs\cO(-2));
 $$
 but this is isomorphic to
 $$
  H^0(X, E_1\otimes E_2^\ast)\otimes H^1(\dP^1,\cO(-2))
 \cong H^0(X,E_1\otimes E_2^\ast),
 $$
  by means of  the K\"{u}nneth formula, the fact that
 $H^0(\dP^1,\cO(-2))=0$ and
 $$H^1(\dP^1,\cO(-2))\cong H^0(\dP^1,\cO)^\ast\cong \dC\ .$$
 Therefore after fixing an element in $H^1(\dP^1,\cO(-2))$, the
 homomorphism $\Phi$ can be identified with the extension class defining $F$.

 Certainly, since the action of $\su$ on the extension class is trivial (note
that this
 action is induced
 from the action on $E_1\otimes E_2^\ast$, which is trivial), then the bundle
$F$ defined
 by the triple $\tri$ is an \sue\ \hvb.
 One can show with a little bit more work that in fact  every \sue\ holomorphic
structure
 on the \sue\ $\ci$ bundle (\ref{equi-f}) defines an extension of the form
 (\ref{extension}).
 \qed

 What this Proposition says is that holomorphic triples over $X$ can be
regarded
 as a ``{\em dimensional reduction}'' of certain \sue\ \hvbs\ over $\xp$.

\subsection{Subtriples}

 \begin{definition} A triple $T' = \stri$\ is a subtriple
 of $T =\tri$\ if
 \begin{tabbing}
 \ \ \ \=(1)\ \= \kill
 \> {\em (1)} \> $E'_i$\ is a coherent subsheaf of $E_i$, for $i=1,2$\\
 \> {\em (2)}\> $\Phi' = \Phi|_{E_2}$, i.e. $\Phi'$\ is the restriction of
 $\Phi$.
 \end{tabbing}
 In other words we have the commutative diagram
 $$
 \begin{array}{ccc}
 E_2& \stackrel{\Phi}{\lra} &  E_1 \\
 \uparrow&  &  \uparrow\\
 E_2'& \stackrel{\Phi'}{\lra} &  E_1'.
 \end{array}
 $$

 If $E'_1=E'_2=0$, the subtriple is called the trivial subtriple.
 \end{definition}

 \noindent {\em Remark.} When studying stability criteria, it will suffice,
  as usual,
 to consider saturated subsheaves, that is subsheaves whose quotient sheaves
are torsion
 free. On a Riemann surface these are precisely subbundles.

 With this definition, subobjects of the triple $\tri$\ are related
 to  subsheaves of the corresponding $SU(2)$-equivariant bundle
 $F\lra \xp$\ in an appropriate way.
 First note that the correspondence between triples on $X$ and bundles on $\xp$
can be
 extended more generally to arbitrary coherent sheaves. Namely, if
 $S_1$ and $S_2$ are two coherent sheaves on $X$ and $\Psi\in \Hom( S_2,S_1)$
the triple
 $( S_1, S_2,\Psi)$
 defines a coherent sheaf $U$ over $\xp$. This sheaf is given, as for
 bundles, as an extension
 \be
 0\lra \ps S_1\lra U\lra\ps S_2\otimes\cod\lra 0.
 \label{c-ext}
 \ee
 The proof is very much as for the case of bundles, once we have fixed
 $ S_1$ and $ S_2$, the extensions as
 (\ref{c-ext}) are parametrized by
 $$
 \ext^1_\xp(\ps S_2\otimes\qs\cod,\ps S_1).
 $$
 But, by the K\"unneth formula for the $\ext$ groups, this group is
 isomorphic to
 $$
 \Hom_X( S_2, S_1)\otimes\ext^1_{\dP^1}(\cO,\cod)\oplus
 \ext^1_X( S_2, S_1)\otimes\Hom_{\dP^1}(\cod,\cO).
 $$
 This reduces to $\Hom_X(S_2,S_1)$, since $\Hom_{\dP^1}(\cod,\cO)\cong
H^0(\dP^1,\cO(-2))=0$ and
 $$\ext^1_{\dP^1}(\cO,\cod)\cong H^1(\dP^1,\cO(-2))\cong\dC\ .$$

 \begin{lemma}\label{equi-sub}
 Let $F\lra\xp$ be the bundle  associated to a triple $\tri$. Then every \sui\
 coherent subsheaf $F'\subset F$ is an extension of the form
 \be
 \sextn,\label{s-ext}
 \ee
 with  $E_1'\subset E_1$ and $E_2'\subset E_2$ coherent subsheaves, making the
 following diagram commutative
 $$
 \begin{array}{ccccccccc}
 0& \lra & \ps E_1 &\lra &  F & \lra & \ps  E_2\otimes\qs\cO(2) &\lra
 & 0\\
  &      &\uparrow  &     & \uparrow & & \uparrow                  &
 &   \\
 0& \lra & \ps E_1' &\lra &  F' & \lra & \ps  E_2'\otimes\qs\cO(2)
 &\lra & 0.
 \end{array}
 $$
  Thus  $F'$ corresponds to a triple $\stri$, for $\Phi'\in \Hom(E_2',E_1')$.
 \end{lemma}
 \pf
 Let $f: F'\ra\ps E_2\otimes\qs\cod$ be the composition of the
 injection $ F'\ra F$ with the surjective map
 $ F\ra\ps E_2\otimes\qs\cod$. Consider the commutative  diagram
 $$
 \begin{array}{ccccccccc}
 0& \lra & \ps E_1 &\lra &  F & \lra & \ps  E_2\otimes\qs\cod &\lra &
 0\\
  &      &\uparrow  &     & \uparrow & & \uparrow                  &
 &   \\
 0& \lra & \Ker f &\lra &  F' & \lra & \im f &\lra & 0.

 \end{array}
 $$
 The $\su$-invariance of $ F'$ implies that of $\Ker f$ and $\im f$.
 It suffices therefore  to show that if $E$ is  a holomorphic vector bundle
over $X$ and
 if $\ps E$ is  the
 pull-back to $\xp$, then every \sui\
 subsheaf of $\ps E$ is isomorphic to a sheaf of the form $\ps E'$ for
 $ E'$ a subsheaf of $ E$.
 Indeed, the action of $\su$ on $\ps E$ can be extended to an action of
$SL(2,\dC)$.
 Let $F'\subset \ps E$ be a $SL(2,\dC)$-invariant coherent subsheaf.
 Consider the action of a subgroup $\dC^\ast\subset SL(2,\dC)$ on
 $X\times\dC\subset\xp$ and  let $A=H^0(X,E)$\ be the space of global sections.
Clearly $H^0(X\times \dC,F')\subset A[t]$, that
 is, $H^0(X\times\dC, F')=\bigoplus_{k=0}^N B_k$, where an element of $B_k$ is
of the
 form $st^k$ for $s\in A$. The action of $\alpha\in \dC^\ast$ is given by
 $$
 \alpha(st^k)=s\alpha^kt^k.
 $$
 By choosing another subgroup $\dC^\ast\subset SL(2,\dC)$, the
$SL(2,\dC)$-invariance of
 $F'$ implies that $H^0(X\times\dC, F')=B_0$ and hence $F'=\ps E'$ for
$E'\subset E$
 a coherent subsheaf.
 \qed

 We shall show in the next lemma that  the triple associated to $F'\subset F$
is in fact
 a subtriple of $\tri$, and conversely, every subtriple of $\tri$ defines a
unique
 \sui\ coherent subsheaf of $F$.

 \begin{lemma}\label{push-out}
 Let $ E_1'\subset E_1$ and $ E_2'\subset E_2$ be coherent subsheaves
 and let $\Phi'\in \Hom( E_2', E_1')$.
 Let $ F'$ be the coherent sheaf over $\xp$ defined by the triple
 $\stri$.
 Then $ F'$ is a subsheaf of $F$ making the  diagram
 $$
 \begin{array}{ccccccccc}
 0& \lra & \ps E_1 &\lra &  F & \lra & \ps  E_2\otimes\qs\cO(2) &\lra
 & 0\\
  &      &\uparrow  &     & \uparrow & & \uparrow                  &
 &   \\
 0& \lra & \ps E_1' &\lra &  F' & \lra & \ps  E_2'\otimes\qs\cO(2)
 &\lra & 0,
 \end{array}
 $$
 commutative if and only if $\stri$ is a subtriple of $\tri$.
 \end{lemma}
 \pf
 Consider the diagram
 $$
 \Hom( E_2', E_1')  \stackrel{i}{\lra}\Hom( E_2',
 E_1)\stackrel{j}{\longleftarrow} \Hom( E_2, E_1).
 $$
 To say that $\stri$ is a subtriple of $\tri$ is equivalent to saying
 that
 $$
 i(\Phi')=j(\Phi).
 $$
 Under the isomorphisms
 $$
 \begin{array}{lll}
 \Hom( E_2', E_1')&\cong &\ext^1(\ps E_1',\ps E_2'\otimes\qs\cod)\\
 \Hom( E_2', E_1)&\cong &\ext^1(\ps E_1,\ps E_2'\otimes\qs\cod)\\
 \Hom( E_2, E_1)&\cong &\ext^1(\ps E_1,\ps E_2\otimes\qs\cod),
 \end{array}
 $$
 $i(\Phi')$ defines an extension $\tilde F^{(i)}$ which makes  the following
 diagram commutative
 \be
 \begin{array}{ccccccccc}
 0& \lra & \ps E_1 &\lra &\tilde F^{(i)} & \lra & \ps  E_2'\otimes\qs\cO(2)
 &\lra & 0\\
  &      &\uparrow  &     & \uparrow & & \parallel                  &
 &   \\
 0& \lra & \ps E_1' &\lra &  F' & \lra & \ps  E_2'\otimes\qs\cO(2)
 &\lra & 0.
 \end{array}
 \label{push-out1}
 \ee
 In particular $ F'$ is a subsheaf of $\tilde F^{(i)}$. On the other hand
 $j(\Phi)$ defines an extension
  $\tilde F^{(j)}$ which fits in the following commutative diagram
 \be
 \begin{array}{ccccccccc}
 0& \lra & \ps E_1 &\lra &  F & \lra & \ps  E_2\otimes\qs\cO(2) &\lra
 & 0\\
  &      &\parallel  &     & \uparrow & & \uparrow                  &
 &   \\
 0& \lra & \ps E_1 &\lra & \tilde F^{(j)} & \lra & \ps
 E_2'\otimes\qs\cO(2) &\lra & 0,
 \end{array}
 \label{push-out2}
 \ee
 and in  particular $\tilde F^{(j)}$ is a subsheaf of $ F$.
 Since $i(\Phi')=j(\Phi)$, $\tilde F^{(i)}\cong \tilde F^{(j)}$ and we can
 compose the above two diagrams to obtain the desired result.
 \qed

 \subsection{Simple Triples}

 \begin{definition}\label{def-simple}Let
 \be
 H^0(E_1,E_2,\Phi)=\bigl\{(u,v)\in H^0(\End E_1)\oplus
 H^0(\End E_2)\bigm|u\Phi=\Phi v\bigr\}.\label{}
 \ee
 We say a holomorphic triple $(E_1,E_2,\Phi)$\ is simple if
 $H^0(E_1,E_2,\Phi)\simeq\Bbb C$, i.e. if the only elements in
 $H^0(E_1,E_2,\Phi)$\ are of the form $\lambda(\bold I_1,\bold I_2)$\
 where $\lambda$\ is a constant and $(\bold I_1,\bold I_2)$\ denote the
identity
 maps on $E_1$\ and $E_2$.
 \end{definition}

 This definition too is dictated  by the correspondence between
 triples on $X$\ and equivariant holomorphic extensions over
 $\XxP$:

 \begin{prop}
 The triple $\tri$ is simple if and only if the \sue\ bundle $F$ associated to
 $\tri$ is $\su$-equivariantly simple.
 \end{prop}
 \pf
 The definition of simplicity for an equivariant bundle is the obvious
 generalization of that for an ordinary holomorphic bundle. Namely, an
equivariant
 bundle is said to be equivariantly simple if it has no other invariant
endomorphisms
 than the constant multiples of the identity.
 The proof of the Proposition follows from the following lemma.
 \begin{lemma}\label{equi-hom}
 Let $T=\tri$ and $T'=\stri$ be two holomorphic triples over $X$ and $F$ and
$F'$ be
 the corresponding \sue\ holomorphic vector bundles over $\xp$. Then, every
\sue\ homomorphism  $g: F\lra F'$
 induces homomorphisms $u: E_1\lra E_1'$ and $v:E_2\lra E_2'$ such that
 \be
 u\Phi=\Phi' v. \label{commu}
 \ee
 Conversely, given morphisms $u$ and $v$ satisfying (\ref{commu}) there exists
a unique
 morphism $g: F\lra F'$ inducing $u$ and $v$.
 \end{lemma}

 \pf
 The map $g$ can be decomposed as
 $$
 g=\left( \begin{array}{cc}
 g_1 & f_1\\
 f_2 & g_2
 \end{array}
 \right),
 $$
 where $g_1:\ps E_1\lra\ps E_1'$, $g_2:\ps E_2\otimes\qs\cod\lra \ps
E_2'\otimes\qs\cod$,
 $f_1:\ps E_2\otimes\qs\cod\lra\ps E_1'$ and $f_2:\ps E_1\lra \ps
E_2'\otimes\qs\cod$.
 By invariance it is very easy to see (cf. \cite[Proposition 3.9]{GP3})  that
 $g_1=\ps u$,\ $g_2=\ps v$, and $f_1=0=f_2$. Equation (\ref{commu}) follows
from
 the commutativity of the following diagram
 $$
 \begin{array}{ccccccccc}
 0& \lra & \;\ps E_1 &\lra &  \;F & \lra & \ps  E_2\otimes\qs\cO(2) &\lra
 & 0\\
  &      &\ps u{\downarrow}  &     & g{\downarrow} & & \ps v{\downarrow}
          &
 &   \\
 0& \lra & \ps E'_1 &\lra &  F' & \lra & \ps E_2'\otimes\qs\cO(2) &\lra & 0.
 \end{array}
 $$
 \qed

 From this  lemma it follows that every \sui\ endomorphism
 of $F$ induces endomorphisms $u: E_1\lra E_1$ and $v: E_2\lra E_2$, satisfying
 $u\Phi=\Phi v$. And conversely, given endomorphisms $u$ and $v$ satisfying
$u\Phi=\Phi v$
 there exists  a unique endomorphism of $F$ inducing $u$ and $v$.
 \qed

 The definition of simplicity for a triple given above  is also motivated by a
deformation
 theory  description of the ``tangent space'' to the space of
 triples.  We will say  more about this in Section \ref{moduli}.

\subsection{Reducible Triples}

 A related, but inequivalent, notion to simplicity is that of
 irreducibility.  We make the following definitions.

 \begin{definition} We say the triple $T=\tri$\ is reducible if
 there are direct sum decompositions $E_1 = \bigoplus_{i=1}^n
 E_{1i}$, $E_2 = \bigoplus_{i=1}^n E_{2 i}$, and
 $\Phi =\bigoplus_{i=1}^n \Phi_i$, such that $\Phi_i\in \Hom(E_{2i},E_{1i})$.
 We adopt the convention that if $E_{2i}=0$\ or
 $E_{1i}=0$\ for some i, then $\Phi_i$\ is the zero map.  With
 $T_i=(E_{1i},E_{2i},\Phi_i)$, we write $T=\bigoplus_{i=1}^n T_i$.
 Thus $T$\ is reducible if it has a decomposition as a direct
 sum of subtriples.

 If $T$ is not reducible, we say $T$ is irreducible.
 \end{definition}

 \begin{prop}
 If a triple $T=\tri$\ is simple, then it is irreducible.
 \end{prop}

 \pf Suppose $T$\ is reducible, with $T=\bigoplus_{i=1}^n T_i$.
 Then we can define $(u,v)\in H^0(E_1,E_2,\Phi)$\ by
 $u=\bigoplus_{i=1}^n \lambda_i \bold I_{1i}$,
 $v=\bigoplus_{i=1}^n \lambda_i  \bold I_{2i}$, where for each $i$,
 $\lambda_i\in \Bbb C$\
 and $\bold I_{1i}(\bold I_{2i})$\ is the identity map on
 $E_{1i}(E_{2i})$. Clearly $T$\ is not simple.
 \qed

 \begin{prop} A holomorphic triple $T=\tri$\ over $X$\ is
 irreducible if and only if the corresponding SU(2)-equivariant
 extension  $F\lra\XxP$\ is equivariantly irreducible, i.e. cannot be
 decomposed as a sum of SU(2)-equivariant extensions of the form (\ref{s-ext}).
 \end{prop}
 \pf This follows directly from the relation between
 subtriples of $T$\ and  \linebreak
$SU(2)$-equivariant subbundles of $F$
 (cf. Lemmas \ref{equi-sub} and \ref{push-out}).
 \qed

 \subsection{Equations for special Metrics}

 Given a \hvb\ over a \crs\ there is a natural condition for a \hm\ on it: that
 of being projectively flat. By choosing a metric on $X$ one can rewrite this
condition
 in a way which turns out  to be the right generalization for higher
dimensional
 manifolds: the {\em \he} condition. Since we shall use this notion on $X$ as
well as on
 $\xp$, we shall define it on a compact \kahler\ manifold of arbitrary
dimension $(M,\omega)$.

 Let $E$ be a \hvb\ over $M$ and $h$ be a \hm\ on $E$. Recall that there is on
$E$ a unique connection
 compatible with both the metric and the holomorphic structure---the so-called
{\em metric
 connection}. Let $F_h$ be its curvature  and $\Lambda F_h$ be the
 contraction of $F_h$ with the \kahler\ form $\omega$. $\Lambda F_h$ is hence a
smooth section
 of $\End E$. The metric $h$ is said to be \he\ with respect to $\omega$ if
 \be
 \sqrt{-1}\Lambda F_h=\lambda \bold I_E,
 \label{he}
 \ee
 where $\bold I_E\in \Omega^0(\End E)$ is the identity and $\lambda$ is a
constant which is
  determined by integrating the trace of (\ref{he}).
 Using  that the degree of $E$, defined as
 $$
 \deg E=\frac{1}{(m-1)!}\int_M c_1(E)\wedge\omega^{m-1},
 $$
 where $m$ is the dimension of  $M$ and $c_1(E)$ is the first Chern class of
$E$,
  given via Chern--Weil theory by
 $$
 \deg E=\frac{i}{2\pi}\int_M\tr(\Lambda F_h) \frac{\omega^m}{m!},
 $$
  we obtain

 $$
 \lambda=\frac{2\pi}{\vol M}\frac{\deg E}{\rank E}.
 $$

 Coming  back to our \crs\ $X$,  let us choose a metric on $X$ with \kahler\
form
  $\omega_X$ and volume normalized to one. Given a \HT\ $\tri$ on $X$ it was
shown in
\cite{GP3} that there are
 natural equations for metrics on the bundles $E_1$ and $E_2$. These equations,
formally
 similar to the \he\ equations,  involve in a natural way  the endomorphism
$\Phi$.
  If $E_1$ and $E_2$ are endowed with Hermitian metrics one can
 form smooth sections of \ $\End E_1$\ and \ $\End E_2$ \ respectively
 by taking the compositions  $\Phi\Phi^\ast$\ and $\Phi^\ast\Phi$.  Here
$\Phi^\ast$\ is the adjoint of $\Phi$\ with respect to the metrics of $E_1$ and
$E_2$.  The equations for the metrics $h_1$ and $h_2$ on $E_1$ and $E_2$,
respectively, are given by

   \be
   \left. \begin{array}{l}
  \sqrt{-1}  \Lambda F_{h_1}+\Phi\Phi^\ast=2\pi\tau \bold I_{E_1}\\
  \sqrt{-1} \Lambda F_{h_2}-\Phi^\ast\Phi=2\pi\tau'\bold I_{E_2}
   \end{array}\right \}, \label{cves}
   \ee
where $\tau$ and $\tau'$ are real parameters.

    We  first observe that, in order to solve (\ref{cves}), the
   parameters $\tau$ and $\tau'$ must be related. Indeed, by adding the trace
of
   the  two equations in  (\ref{cves}), and since $\tr(\Phi\Phi^\ast)=
   \tr(\Phi^\ast\Phi)$, we get
   $$
   \sqrt{-1}\tr(\Lambda F_{h_1})+\sqrt{-1} \tr(\Lambda F_{h_2})=2\pi r_1\tau
+2\pi r_2\tau',
   $$
 where $r_1$ and $r_2$ are the ranks of $E_1$ and $E_2$ respectively.
   By integrating this equation  and, since
   $$
   \deg E_1=\frac{\sqrt{-1}}{2\pi}\int_X \tr(\Lambda F_{h_1})\omega
   \;\;\;\mbox{and}\;\;\;
   \deg E_2=\frac{\sqrt{-1}}{2\pi}\int_X \tr(\Lambda F_{h_2})\omega,
   $$
   we obtain
 \be
   r_1\tau+r_2\tau '=\deg E_1+\deg E_2.
   \label{t-t'}
 \ee

 There is therefore just one independent  parameter, that we choose to be
$\tau$.
 These equations are called the {\em coupled $\tau$-vortex equations}
 by analogy with the vortex equations on a single bundle studied in
\cite{B1,B2,GP1,GP2}.

\subsection{Dimensional Reduction of Equations }

 The coupled vortex equations have similar interpretations to the \he\ equation
 and the vortex equations on a single bundle. They can be interpreted both as
 the equations satisfied by the minima of a certain gauge-theoretical
 functional---a generalized \ymh-type
 functional---as well as moment map equations in the sense of symplectic
 geometry (see \cite[Section 2]{GP3} for details).
 In fact, the relation between  the \cves\ and the \he\ equation that we shall
exploit here is of a more intimate
 nature. Namely, the \cves\ are a dimensional reduction of the \he\ equation
under the
 action of $\su$ on $\xp$.
 Of course, in order to talk about the \he\ equation on $\xp$ one needs to
choose a \kahler\ metric.
 We shall consider the one-parameter family  of  \sui\ \kahler\ metrics
 with \kahler\ form
 $$
 \os=\frac{\sigma}{2}\ps\omega_X\oplus\omega_{\dP^1},
 $$
 where $\omega_{\dP^1}$ is the Fubini-Study \kahler\ form normalized to volume
one, and $\sigma\in\dR^+$.

 \begin{prop}\label{dr}Let $T=\tri$ be a holomorphic triple and $F$ be the
\sue\
 holomorphic bundle over $\xp$ associated to $T$, that is given as an extension
 \be
 \extn.\label{extension-a}
 \ee
 Suppose  that $\tau$ and $\tau'$ are related by (\ref{t-t'}) and  let
 \be
 \sig=\frac{(r_1+r_2)\tau-(\deg E_1 +\deg  E_2)}{r_2}.\label{s-t}
 \ee
 Then $E_1$ and $E_2$ admit metrics satisfying the \ctves\ if and
 only if $F$ admits an \sui\ \he\ metric with respect to $\os$.
 \end{prop}
 \pf
 We shall give here just a sketch of the proof (see \cite[Proposition
3.11]{GP3} for details).
 First one has the following result, which is a special case of the general
characterization of an \sui\ Hermitian metric on an  \sue\ vector bundle over
$\xp$.

 \begin{lemma} Let $h$ be an \sui\ \hm\ on the bundle $F\lra\xp$ associated to
the triple
 $\tri$. Then $h$ is of the form
 \be
 h=\ps h_1\oplus\ps h_2\otimes\qs h_2',\label{inv-met}
 \ee
 where $h_1$ and $h_2$ are metrics on $E_1$ and $E_2$, respectively, and
 $h_2'$ is an \sui\ metric on $\cod$.
 Conversely, given metrics $h_1$, $h_2$ and $h_2'$ as above, (\ref{inv-met})
defines
 an \sui\ metric on $F$.
 \end{lemma}

 Let $F_1$ and $F_2$ be the curvatures of the metric connections of
 $\ps h_1$ and $\ps h_2\otimes \qs h_2'$ respectively. Then
 $$
 \begin{array}{l}
 F_1=\ps F_{h_1}\\
 F_2=\ps F_{h_2}\otimes 1 + \bold I_{E_2}\otimes\qs F_{h_2'}.
 \end{array}
 $$
 The curvature of the metric connection corresponding to $h$ is given by
  $$
    F_h= \left( \begin{array}{cc}   F_1 - \beta
 \wedge \beta ^\ast & D' \beta \\
 -D''\beta^\ast & F_2 -\beta ^\ast
   \wedge \beta \end{array} \right),
  $$
 where $\beta\in\Omega^{0,1}(\xp,\ps(E_1\otimes E_2^\ast)\otimes\qs\cO(-2))$ is
a
 representative of the extension class in $H^1(\xp,\ps(E_1\otimes
E_2^\ast)\otimes\qs\cO(-2))$
 defining (\ref{extension-a}), and
 $$
 D:\Omega^1(\xp,\ps(E_1\otimes E_2^\ast)\otimes\qs\cO(-2))\lra
   \Omega^2(\xp,\ps(E_1\otimes E_2^\ast)\otimes\qs\cO(-2))
 $$
 is built from the metric connections of $\ps h_1$ and $\ps h_2\otimes \qs
h_2'$.

 As explained in Proposition \ref{equi-hvb}, $\beta=\ps\Phi\otimes\qs\alpha$,
where
 $\alpha\in\Omega^{0,1}(\dP^1,\cO(-2))$ is the unique \sui\ representative of
the
 element in $H^1(\dP^1,\cO(-2))$, which has to be fixed in order to associate
the
 extension (\ref{extension-a}) to $\tri$. One can choose the constant in
 $H^1(\dP^1,\cO(-2))\cong\dC$ such that
$\alpha\wedge\alpha^\ast=\frac{1}{\sigma}\omega_{\dP^1}$.

 Let $\Lambda_\sigma$ be the contraction with the \kahler\ form $\os$. A
straightforward
 computation shows that if $\sigma$ is related to $\tau$ by (\ref{s-t}), then
$h$
 is \he\ with respect to $\os$. That is
 $$
 \sqrt{-1}\Lambda_\sigma F_h=\lambda\bold I_F,
 $$
 if and only if $h_1$ and $h_2$ satisfy the \ctves.

 We have assumed that if  the relation between  $\sigma$ and  $\tau$ is  given
by
 (\ref{s-t}), then $\sigma>0$. However, we will show
 in Section  \ref{stability} that this can actually be derived from the \cves.
 \qed

 \noindent{\em Remark}.\ The choice of the \kahler\ metric on $\xp$ that we
have made  differs
 from the one made in \cite{GP3}. There the parameter $\sigma$ is multiplying
the metric on $\dP^1$,
  i.e.
 $\os=\ps\omega_X\oplus\sigma\qs\omega_{\dP^1}$. This, and the fact that the
volume of $X$
 was not  normalized to one, explains why the relation between $\tau$ and
$\sigma$ given there is  the inverse of (\ref{s-t}).

\subsection{Invariant Stability and the Hitchin-Kobayashi
Correspondence}

It is very well-known that the existence of a \he\ metric on a \hvb\ is
governed by
 the algebraic-geometric condition of {\em stability}. Recall that a \hvb\ $E$
over
 a compact \kahler\ manifold $(M,\omega)$ is said to be {\em stable} if
 $$
 \mu(E')<\mu(E)
 $$
 for every non-trivial coherent subsheaf $E'\subset E$. Where
 $$
 \mu(E')=\frac{\deg E'}{\rank E'}
 $$
 is the {\em slope} of $E'$.

 The precise relation between the \he\ condition and stability is given by the
 so-called {\em Hitchin--Kobayashi correspondence}, proved by Donaldson
\cite{D1,D2} in the
 algebraic case and by Uhlenbeck and Yau \cite{U-Y} for an arbitrary compact
\kahler\ manifold (see also \cite{Ko,L,A-B,N-S}):

 \begin{thm}\label{h-k}
 Let $E$ be a \hvb\ over a compact \kahler\ manifold $(M,\omega)$. Then $E$
admits
 a \he\ metric if and only if $E$ is polystable, that is a direct sum of stable
bundles
 of the same slope.
 \end{thm}

 From this theorem and Proposition \ref{dr} we conclude that the existence of
solutions
 to the \cves\ must be dictated by the stability of the bundle $F\lra\xp$
associated to
 the triple $\tri$. In fact, since the \he\  metric on $F$ is \sui, the
condition that $F$
 has to satisfy is a slightly weaker condition than stability,
 namely  that of {\em invariant stability}.
 Let $(M,\omega)$ be a compact \kahler\ manifold and $G$ be a compact Lie group
acting on
 $M$ by isometric biholomorphisms. Let $E$ be a $G$-equivariant \hvb\ over $M$.
We
 say that $E$ is $G$-{\em invariantly stable} if
 $$
 \mu(E')<\mu(E)
 $$
 for every $G$-invariant non-trivial coherent subsheaf $E'\subset E$.

 The basic relation between $G$-invariant stability and ordinary stability is
 given by the following theorem (cf. \cite[Theorem 4]{GP2}).

 \begin{thm}
 \label{isvs}
 Let $E$ be a $G$-invariant holomorphic vector
 bundle as above. Then $E$ is $G$-invariantly stable if and only if
 $E$ is $G$-indecomposable and is of the form
 $$
 E=\bigoplus _{i=1}^n E_i
 $$
 \noindent where $E_i$ is a  stable bundle, which is  the transformed of $E_1$
 by an element of $G$.
 \end{thm}

 As a corollary of Theorems \ref{h-k} and \ref{isvs} one obtains a
$G$-invariant version
 of the Hitchin-Kobayashi correspondence (cf. \cite[Theorems 4 and 5]{GP2}):

 \begin{thm}\label{inv-h-k}
 Let $E$ be a $G$-equivariant  \hvb\ over a compact \kahler\ manifold
$(M,\omega)$. Then $E$
 admits a $G$-invariant \he\ metric if and only if $E$ is $G$-invariantly
polystable, that
 is a direct sum of $G$-invariantly stable bundles of the same slope.
 \end{thm}

 From Proposition \ref{dr} and Theorem \ref{inv-h-k} we obtain the following
existence
 theorem.

 \begin{thm}\label{exst}
Let $T=\tri$ be a \HT\ over a \crs\ $X$ equipped with a metric. Let $F\lra\xp$
 be the bundle associated to $T$ as above. Let $\sigma$ and $\tau$\ be real
parameters related by (\ref{s-t}). Then $E_1$ and $E_2$ admit metrics
satisfying
 the \ctves\ if and only if $F$ is a $\su$-invariantly polystable bundle with
respect
 to the \kahler\ form $\os$ defined above.
 \end{thm}

 \section{Definition and Properties of Stability for Triples}\label{stability}

The existence theorem \ref{exst}  gives conditions on the extension $F\lra\xp$\
for existence of solutions to the coupled vortex equations on $\tri$. We would
like to express these conditions entirely in terms of the data on $\tri$.
Indeed this is one of our primary objectives in this paper. To achieve this, we
will need an appropriate notion of stability for a triple. In this section we
define such a concept  for holomorphic triples, and discuss some properties
that follow from the definition.

 Keeping our earlier notation, we let  $E_1$\ and $E_2$\ be holomorphic vector
bundles over a
  Riemann surface $X$.  We denote their ranks by $r_1$\ and $r_2$\
respectively, and their degrees by $d_1$\
 and $d_2$.  We let $\Phi:E_2\lra E_1$\ be a holomorphic bundle homomorphism,
i.e.
 $\Phi\in H^0(\Hom(E_2,E_1))$.

 Our definition of stability for the triple $T=\tri$\ has two equivalent
formulations.  The first has some advantages when considering the relation
between stability and the coupled vortex equations, while the second  has the
virtue that it is in the style of the definition of parabolic stability, and
thus looks more familiar.  Both definitions involve a real parameter, with the
result that there is a 1-parameter family of stability criteria for triples.
This is the same phenomenon as is observed in the case of holomorphic pairs.
All our results can be of course be
 stated in terms of either definition, and for the sake of completeness we will
give both versions.

 \begin{definition}\label{t-stab} Let \ $T'=(E'_1,E'_2,\Phi)$\ be a nontrivial
subtriple of
\ $(E_1,E_2,\Phi)$, with \ $\rank E'_1=r_1'$\ and \ $\rank E'_2=r_2'$.  For any
real $\tau$\ define

 \be
 \theta_{\tau}(T')=(\mu(E'_1\oplus
E'_2)-\tau)-\frac{r_2'}{r_2}\frac{r_1+r_2}{r_1'+r_2'}
 (\mu(E_1\oplus E_2)-\tau).\label{Theta}
 \ee
 The triple $T=(E_1,E_2,\Phi)$ is called $\tau$-stable if
 $$\theta_{\tau}(T')<0$$
 for all nontrivial subtriples $T'=(E'_1,E'_2,\Phi)$.
 The triple is called $\tau$-semistable if for all subtriples
 $$\theta_{\tau}(T')\le0.$$
 \end{definition}

 \begin{definition}\label{s-stab} With $\sigma$\ a real number, define the
 $\sigma$-degree and $\sigma$-slope of a subtriple $T'=(E'_1,E'_2,\Phi)$\ by
 $$\deg_{\sigma}(T')=\deg(E_1'\oplus E_2')+r_2'\sigma ,$$
 and
 $$\mu_{\sigma}(T')=\frac{\deg_{\sigma}(T')}{r_1'+r_2'}\ .$$

  The triple $T=(E_1,E_2,\Phi)$\ is called
 $\sigma$-stable if for all nontrivial subtriples
 $T'=(E'_1,E'_2,\Phi)$\ we have
 $$\mu_{\sigma}(T') < \mu_{\sigma}(T)\ .$$
 \end{definition}

 A straightforward computation shows the equivalence of these two definitions.
 \begin{prop}\label{tstab-sstab} Fix $\tau$\ and $\sigma$\ such that
 $$\sigma = \frac{r_1+r_2}{r_2}(\tau-\mu(T))\ ,$$
 or equivalently
 $$\tau = \mu_{\sigma}(T)\ .$$
 Then for any subtriple $T'=(E'_1,E'_2,\Phi)$, the following are
 equivalent:
 \begin{tabbing}
 \ \ \ \=(1)\ \= \kill
 \> {\em (1)}\> $\theta_{\tau}(T')<0\ $,\\
 \> {\em (2)} $\mu_{\sigma}(T') < \mu_{\sigma}(T)$.
 \end{tabbing}
 That is, the triple is $\tau$-stable if and only if it is
 $\sigma$-stable.
 A similar result holds with ``~ $<$\ ~" replaced by ``~ $=$\ ~".
 \end{prop}

 \noindent {\em Remark.}  There are two special cases where the notion of
stability for a triple is especially simple, namely when $\Phi=0$, and when
$E_2$\ is a line bundle.

 \begin{lemma}Suppose that $\Phi=0$.  The degenerate holomorphic triple
$(E_1,E_2,0)$\ is $\tau$-semistable if and only if $\tau=\mu(E_1)$\ and both
bundles are semistable.  Such triple cannot be $\tau$-stable.
 \end{lemma}
 \pf Subtriples of $T=(E_1,E_2,0)$\ are all of the form $T'=(E'_1,E'_2,0)$,
with $E'_1$\ and $E'_2$\ being
 any holomorphic subbundles of $E_1$\ and $E_2$\ respectively.  Applying the
condition $\theta_{\tau}(T')\le 0$\ to subtriples of the form $T'=(E'_1,0,0)$\
gives
 \be
 \mu(E'_1)\le \tau,\label{3.2}
 \ee
 while applying the condition to subtriples of the form $T'=(E'_1,E_2,0)$\
gives
stability
 \be
 \mu(E_1/E'_1)\ge \tau.\label{3.3}
 \ee
 These two inequalities imply
 \be
 \mu(E_1)\le\tau\le\mu(E_1).\label{3.4}
 \ee
 That is, $\tau=\mu(E_1)$, and hence $E_1$\ is a semistable bundle.  Similarly,
by considering the subtriples
  $(0,E'_2,0)$\ and $(E_1,E'_2,0)$, we see that $E_2$\ is also semistable.
Notice that the inequalities in
 (\ref{3.2}) and (\ref{3.3}) cannot be made strict without leading to a
contradiction in (\ref{3.4}).
 \qed

 \begin{cor}The map $\Phi$\ cannot be identically zero in a $\tau$-stable
triple.
 \end{cor}

 \begin{lemma}In the case where $E_2=L$\ is a line bundle, i.e. $r_2=1$, the
above definition is
 equivalent to the notion of $\tau$-stability defined in \cite{GP3}. It thus
corresponds to the  $(\tau-\deg L$)-stability for the holomorphic pair
$(E_1\otimes L^*, \Phi)$.
 \end{lemma}
 \pf In this case there are only two types of subtriple possible, corresponding
to $r_2'=0$\ or $r_2'=1$.
 In the first case the subtriples are of the form $(E'_1,0,0)$, where $E'_1$\
is an arbitrary holomorphic subbundle of $E_1$.  The condition
$\theta_{\tau}(T')<0$\
 then reduces to
 $$\mu(E'_1)<\tau.$$
 In the second case, the subtriples are of the form $(E'_1,E_2,\Phi)$\ where
$E'_1$\ is a holomorphic subbundle such that $\Phi(E_2)\subset E'_1$.  For such
subtriples the condition $\theta_{\tau}(T')<0$\ is equivalent to
 $$(r_1'+1)\mu(E'_1\oplus E_2)-(r_1+1)\mu(E_1\oplus E_2)-(r_1'-r_1)\tau <0,$$
 i.e.
 $$\mu(E_1/E'_1)>\tau.$$
 \qed

 Definition \ref{t-stab} can thus be considered a natural extension of the
$\tau$-stability for pairs defined in
 \cite{B2}. For the more general triples which we are considering here however,
the number of different
 possibilities for subtriples is too large to reformulate the definition of
$\tau$-stability in the style of \cite{GP3} or \cite{B2}, i.e. in terms of
separate slope conditions on the various families of subtriples.  The
$\tau$-stability\ of a triple does however imply the following conditions on
subtriples:

 \begin{prop}\label{slopes}Let $(E_1,E_2,\Phi)$\ be a $\tau$-stable triple. Let
$\tau'$\ be related to $\tau$ by
 \be
 r_1\tau+r_2\tau'=\deg E_1+\deg E_2.\label{3.5}
 \ee
 Then
 \begin{tabbing}
 \ \ \ \=(0)\ \= \kill
 \> {\em (1)}\> $\mu(E_1')<\tau$\ for all holomorphic subbundles $E_1'\subset
E_1$,\\
 \> {\em (2)}\>  $\mu(E_2')<\tau'$\ for all holomorphic subbundles $E_2'\subset
E_2$\ such that $E_2'\subset \Ker(\Phi),$\\
 \> {\em (3)}\>  $\mu(E_2'')>\tau'$\ for all holomorphic quotients of $E_2$,\\
 \> {\em (4)}\> $\mu(E_1'')>\tau$\ for all holomorphic quotients of $E_1$\ such
that $\pi\circ\Phi(E_2)=0$,\\
 \> \> where $\pi:E_1\lra E_1''$\ denotes projection onto the quotient.
 \end{tabbing}
 \end{prop}
 \pf
 These are immediate consequences of the stability condition, i.e.
 $$\theta_{\tau}(E'_1,E'_2,\Phi')<0,$$
 applied to the following special subtriples
 \begin{tabbing}
 \ \ \ \=(0)\ \= \kill
 \> (1) $(E_1',0,\Phi)$,\\
 \> (2) $(0,E_2',\Phi)$,\\
 \> (3) $(E_1,E'_2,\Phi)$, with $E_2''=E_2/E'_2$,\\
 \> (4) $(E'_1,E_2,\Phi)$, with $E_1''=E_1/E'_1$.
 \end{tabbing}
 \qed

 Notice that (\ref{3.5}) can be expressed as
 $$\tau' = \mu(E_1\oplus E_2) - \frac{r_1}{r_1+r_2}\sigma\ .$$
 An equivalent formulation of Proposition \ref{slopes} is thus

 \begin{prop}\label{slopes'}Let $T=(E_1,E_2,\Phi)$\ be a $\sigma$-stable
triple.  Then
 \begin{tabbing}
 \ \ \ \=(0)\ \= \kill
 \> {\em (1)}\> $\mu(E_1')<\mu(T) + \frac{r_2}{r_1+r_2}\sigma$\ for all
holomorphic subbundles $E_1'\subset E_1$,\\
 \> {\em (2)}\> $\mu(E_2')<\mu(T) - \frac{r_1}{r_1+r_2}\sigma$\ for all
holomorphic subbundles $E_2'\subset E_2$\
 such that\\
 \> \> $E_2'\subset \Ker\Phi,$\\
 \> {\em (3)}\> $\mu(E_2'')>\mu(T) - \frac{r_1}{r_1+r_2}\sigma$\ for all
holomorphic quotients, $E_2''$, of $E_2$,\\
 \> {\em (4)}\> $\mu(E_1'')>\mu(T) + \frac{r_2}{r_1+r_2}\sigma$\ for all
holomorphic quotients, $E_1''$, of $E_1$\
 such that\\
 \> \> $\pi\circ\Phi(E_2)=0$, where $\pi:E_1\lra E_1''$\ denotes projection
onto the quotient.
 \end{tabbing}
 \end{prop}

 \subsection{Stable implies simple}

 An important consequence of stability for holomorphic bundles is that the only
automorphisms of a stable bundle are the constant multiples of the identity,
i.e. stable bundles are simple.  We now show that this remains true in the case
of holomorphic triples, where the definition of simplicity is that given in
Definition \ref{def-simple}.  The key result is the following Proposition.

 \begin {prop}\label{s-s} Let $(E_1,E_2,\Phi)$\ be a $\tau$-stable holomorphic
triple. Let $(u,v)$\ be in $H^0(E_1,E_2,\Phi)$.  Either  $(u,v)$\ is trivial,
or both $u$\ and $v$\ are isomorphisms.
 \end{prop}
 \pf Suppose that $u$\ and $v$\ are both neither trivial nor isomorphisms.
Consider the triples
 $K=(\Ker u,\Ker v,\Phi)$\ and $I=(\im u,\im v,\Phi)$, where $\Ker$\ and $\im$\
denotes the kernels and
 images of the maps.  Since $u\Phi=\Phi v$,  these are both proper subtriples
of
 $\tri$, and thus the $\tau$-stability condition gives
 \be
 \theta_{\tau}(K)<0,\label{3.17}
 \ee
 and
 \be\theta_{\tau}(I)<0.\label{3.18}
 \ee
 We also have the exact sequences
 $$
 0\lra\Ker u\lra E_1\lra\im u\lra0,
 $$
 and
 $$
 0\lra\Ker v\lra E_2\lra\im v\lra0.
 $$
 Let $\sigma_u$\ and $\rho_u$\ denote the ranks of $\Ker u$\ and $\im u$, and
similarly for $\sigma_v$\ and $\rho_v$.  Then from the exact sequences we get
 \be
 (\sigma_u+\sigma_v)\mu(\Ker u\oplus \Ker v) + (\rho_u+\rho_v)\mu(\im u\oplus
\im v)=
 (r_1+r_2)\mu(E_1\oplus E_2).\label{3.19}
 \ee
 But by definition of $\theta_{\tau}(K)$,

 $$
 \begin{array}{ll}
 (\sigma_u+\sigma_v)\mu(\Ker u\oplus \Ker v)&=
 (\sigma_u+\sigma_v)\theta_{\tau}(K)\\
 &+\frac{\sigma_v}{r_2}(r_1+r_2)\mu(E_1\oplus E_2)
 +(\sigma_u+\sigma_v-\frac{\sigma_v}{r_2}(r_1+r_2))\tau,
 \end{array}
 $$
 with a similar expression for $(\rho_u+\rho_v)(\mu(\im u\oplus \im v)$.  Also,
 $\sigma_u+\rho_u=r_1$, and $\sigma_v+\rho_v=r_2$.  Hence from (\ref{3.19}) we
obtain
 $$r_1\theta_{\tau}(K)+r_2\theta_{\tau}(I)=0.$$
 This is incompatible with (\ref{3.17}) and (\ref{3.18}).
 \qed

 \begin{cor}\label{stab-simple} If $(E_1,E_2,\Phi)$\ is $\tau$-stable, then it
is simple.
 \end{cor}
 \pf Let $(u,v)$\ be a nontrivial element in $H^0(E_1,E_2,\Phi)$.  By the above
Proposition, both $u$\ and $v$\ are isomorphisms.  Fix a point $p$\ on the base
of the bundles, and let $\lambda$\ be an eigenvalue of $v:E_2|_p\lra E_2|_p$,
i.e. of $v$\ acting on the fibre over $p$.

 Now define
 $$\hat u =u-\lambda \bold I_1,$$
 $$\hat v =v-\lambda \bold I_2.$$

 Clearly $(\hat u,\hat v)$\ is in $H^0(E_1,E_2,\Phi)$, but since $\hat u$\  is
not an
isomorphism,  it follows from Proposition \ref{s-s} that both are identically
zero, i.e.
 $$(u,v)=\lambda (\bold I_1,\bold I_2).$$
 \qed

 We see, in particular, that stable triples are necessarily irreducible.  For
reducible triples, we can however define a notion of polystability. This will
be useful when we consider the relation between stability and the coupled
vortex equations.

 \begin{definition}Let $T=\tri$\ be a reducible triple, with
$T=\bigoplus_{i=1}^n T_i$.  Suppose that
 in each summand $T_i=(E_{1i},E_{2i},\Phi_i)$, the map $\Phi_i$\ is non-trivial
unless $E_{1i}=0$\ or
 $E_{2i}=0$.  Fix value of $\tau$, and let $\tau'$\ be related to $\tau$\ as in
(\ref{3.5}).  We say
 that $T$\ is $\tau$-polystable if for each summand $T_i$
 \begin{tabbing}
 \ \ \ \= (0)\ \= \kill
 \> {\em (1)}\> if $\Phi_i\ne 0$, then $T_i$ is $\tau$-stable,\\
 \> {\em (2)}\> if $E_{1i}=0$, then $E_{2i}$\ is a stable bundle of slope
$\tau'$,\\
 \> {\em (3)}\> if $E_{2i}=0$, then $E_{1i}$\ is a stable bundle of slope
$\tau$.
 \end{tabbing}
 \end{definition}

 \subsection{Duality for triples}

 Associated to a triple $T=\tri$ there is always a  {\em dual triple}
 $T^*=\dtri$, where $\Phi^*$ is the transpose of $\Phi$, i.e.
 the image of $\Phi$ via the canonical isomorphism
 $$
 \Hom(E_2,E_1)\cong \Hom(E_1^*,E_2^*).
 $$
 It is reasonable that the stability of $T$ should be related  to that of
 $T^*$. More precisely.

 \begin{prop}\label{duality}
 $T=\tri$ is \ts\ if and only if $T^*=\dtri$ is $(-\tau')$-stable, where
 $\tau'$ is related to $\tau$ by (\ref{3.5}).
 Equivalently, $T$\ is $\sigma$-stable if and only if $T^*$\ is
$\sigma$-stable.
 \end{prop}
 \pf
 Let $T'=\stri$ be a subtriple of $T$. This defines  a quotient triple
 $T''=(E_1'',E_2'',\Phi'')$, where $E_1''=E_1/E_1'$, $E_2''=E_2/E_2'$, and
 $\Phi''$ is the morphism induced by $\Phi$.
 $T''^*=(E_2''^*,E_1''^*,\Phi''^*)$ is the desired subtriple of
 $T^*$. Since one has the isomorphism $T\cong T^{**}$ we can conclude
 that there is a one-to-one correspondence between subtriples of $T$ and
subtriples of $T^*$. It is not difficult to verify that $\TT(T')<0$ is
equivalent to $\theta_{-\tau'}(T''^*)<0$.  The equivalence of the
$\sigma$-stability for $T$\ and $T^\ast$\ now follows from the fact that if
$\tau=\mu_{\sigma}(T)$, then
 $$-\tau'=-\mu(E_1\oplus E_2)-\frac{r_1}{r_1+r_2}\sigma = -\mu_{\sigma}(T^*)\
.$$
 \qed

 \subsection{Constraints on the parameters}

 \begin{prop}\label{l-bound}Let $(E_1,E_2,\Phi)$\ be a $\tau$-stable triple,
and let $\tau'$\ be as above.  Then
 \begin{tabbing}
 \ \ \ \= (0)\  \= \kill
 \> {\em (1)}\> $\tau >\mu(E_1)$,\\
 \> {\em (2)}\> $\tau' <\mu(E_2)$, and\\
 \> {\em (3)}\> $\tau -\tau' >0$,
 \end{tabbing}
 Equivalently, if $(E_1,E_2,\Phi)$\ is $\sigma$-stable, then
 \begin{tabbing}
 \ \ \ \= (0)\ \=\kill
 \> {\em (1)}\> $\sigma > \mu(E_1)-\mu(E_2)$,\\
 \>  {\em (2)}\> $\sigma > 0$.
 \end{tabbing}
 \end{prop}
 \pf The first two statement follows from cases (1) and (3) in Proposition
\ref{slopes} with $E_1'=E_1$, and
 $E''_2=E_2$\ respectively.

 To prove the third statement, let $K$\ be the subbundle of $E_2$\ generated by
the kernel of $\Phi$, and let $I$\ be the subbundle of $E_1$\ generated by the
image of $\Phi$.  Since the triple is assumed to be $\tau$-stable, $\Phi$, and
therefore $I$, is non-trivial.  By (1) in Proposition \ref{slopes} we thus have
 \be
 \mu(I)<\tau.\label{3.6}
 \ee
 But we also have $0\lra K\lra E_2\lra I\lra 0$, i.e. $I$\ is a quotient of
$E_2$.  It thus follows from (3) in Proposition \ref{slopes} that
 \be
 \mu(I)>\tau'.\label{3.7}
 \ee

 The bounds on $\sigma$\ can be obtained from those on $\tau$\ by substituting
$\tau =\mu_{\sigma}(T)$, and using the fact that
 $$\sigma =\tau-\tau'$$
 if $\tau'$\ is as above.
 \qed

 Part (1) of this proposition gives the lower bound on the allowed range for
$\tau$.  In almost all cases the rank and degree of $E_1$\ and $E_2$\ also
impose an upper bound on $\tau$. In fact

 \begin{prop}\label{u-bound}Let $\tri$\ be a triple with $r_1\ne r_2$.  If the
triple is \ts\, then
 \be
 \tau<\mu(E_1)+\frac{r_2}{|r_1-r_2|}(\mu(E_1)-\mu(E_2))\label{3.8a}
 \ee
 Equivalently, if the triple is $\sigma$-stable, then
 \be
 \sigma < (1+\frac {r_1+r_2}{|r_1-r_2|})(\mu(E_1)-\mu(E_2)).\label{3.8b}
 \ee
 \end{prop}
 \pf
 Let $K=\Ker \Phi$ and $I=\im \Phi$. Consider the subtriples
 $$
 T_1=(0, K,\Phi)\;\;\;\;\text{and}\;\;\;\;T_2=(I,E_2,\Phi).
 $$
 Since $r_1\ne r_2$, $\Phi$\ cannot be an isomorphism and at least one of these
must be a proper subtriple.
 Let $r_2'=\rank K$, $r_2''=\rank I$, $d_2'=\deg K$ and $d_2''=\deg I$.

 A straightforward computation shows that
 \begin{eqnarray}
 \TT(T_1)<0\Longleftrightarrow& d_2'-r_2'(d_1+d_2)+r_1r_2'\tau<0 \label{3.9}\\
 \TT(T_2)<0\Longleftrightarrow&d_2''-d_1+(r_1-r_2'')\tau<0.\label{3.10}
 \end{eqnarray}

 Adding  (\ref{3.9}) to $r_2$ times (\ref{3.10}), and noting that
$d_2=d_2'+d_2''$ and
 $r_2=r_2'+r_2''$, we get that
 \be
 r_2(d_2-d_1)-r_2'(d_1+d_2)+(r_2(r_1-r_2)+r_2'(r_1+r_2))\tau<0.
 \label{3.11}
 \ee
 On the other hand combining (3) in Proposition \ref{l-bound} and (\ref{3.5})
we obtain
 \be
 d_1+d_2-(r_1+r_2)\tau<0.
 \label{3.12}
 \ee

 Adding (\ref{3.11}) to $r_2'$ times (\ref{3.12}) we get
 $$
 (r_1-r_2)\tau<d_1-d_2.
 $$
 If now  $r_1>r_2$, then we get
 \be
 \tau<\frac{d_1-d_2}{r_1-r_2}\label{3.13}
 \ee
 or equivalently
 $$
 \tau<\mu(E_1)+\frac{r_2}{r_1-r_2}(\mu(E_1)-\mu(E_2)).
 $$

 To obtain the bound in the case $r_1<r_2$, note that by Proposition
\ref{duality}
 the $\tau$-stability of $\tri$ is equivalent to the $(-\tau')$-stability of
the dual triple
 $\dtri$, where  $\tau'$ is given by
 \be
 r_1\tau+r_2\tau'=d_1+d_2.
 \label{3.14}
 \ee
 Hence we can  apply  (\ref{3.13}) to $\dtri$ to get that
 $$
 -\tau'<\frac{d_1-d_2}{r_2-r_1},
 $$
 which together with (\ref{3.14}) leads to
 $$
 \tau<\frac{(2r_2-r_1)d_1-r_1d_2}{(r_2-r_1)r_1},
 $$
  i.e.
 $$
 \tau<\mu(E_1)+\frac{r_2}{r_2-r_1}(\mu(E_1)-\mu(E_2)).
 $$
 \qed

 Combining the lower and upper bounds on $\tau$\ (or $\sigma$) we can deduce

 \begin{cor}  If $\rank E_1$ and $\rank E_2$\ are unequal, then a triple
$\tri$\ cannot be stable unless $\mu(E_2)<\mu(E_1)$.
 \end{cor}

 Furthermore,  by the proof of Proposition \ref{u-bound} we get the following
corollary.

 \begin{cor}\label{Corollary 3.8}Let $\tri$\ be \ts\, and suppose that $r=s$.
If $\Phi$\ is not an isomorphism, then $d_1>d_2$.  In particular, in any \ts\
triple $\tri$, the bundle map $\Phi$\ is an isomorphism if and only if
$r_1=r_2$\ and $d_1=d_2$.
 \end{cor}
 \pf The fact that $d_1>d_2$ follows from the inequality
 $$
 (r_1-r_2)\tau<d_1-d_2\ ,
 $$
 which applies if $\Phi$\ is not an isomorphism.  In particular, if
 $\Phi$\ is not an isomorphism then $d_1\ne d_2$.  Conversely, if $\Phi$\ is an
isomorphism, then clearly $r_1=r_2$\ and $d_1=d_2$.
 \qed

 It is when $\Phi$\ is an isomorphism that the range for $\tau$\ can fail to be
bounded.  For example

 \begin{prop}\label{Proposition 3.6}Suppose that $E_1$\ and $E_2$\ are both
stable bundles of rank
 $r$ and degree $d$, and that $\Phi:E_2\lra E_1$\ is non trivial.  Then for any
$\tau>\mu(E_1)$\ the holomorphic triple $(E_1,E_2,\Phi)$\ is $\tau$-stable.
 \end{prop}
 \pf
 Let $\mu(T)=\mu(E_1\oplus E_2)$, and for a subtriple $T'=(E'_1,E'_2,\Phi)$\
set $\mu(T')=\mu(E'_1\oplus E'_2)$.  Since $E_1$\ and $E_2$\ are stable and of
equal slope, we have $\mu(T')<\mu(T)$\ for all subtriples. Thus
 $$\theta_{\tau}(T')\le (\mu(T)-\tau)\frac{r_1'-r_2'}{r_1'+r_2'}.$$
 Since $\Phi$\ is a nontrivial map between stable bundles of the same rank and
degree, it must be a multiple of the identity (cf. \cite{O-S-S}). In particular
$\Phi$\ is injective and hence $r_1'-r_2'\ge 0$.
 Thus $\theta_{\tau}(T')\le 0$.
 In fact, $\theta_{\tau}(T')<0$\ unless $r_1'=r_2'$.  But in that case, we can
write
 $$\theta_{\tau}(T')=r_1'(\mu(E_1')-\mu(E_1))+r_1'(\mu(E_2')-\mu(E_2)),$$
 which is strictly negative.
 \qed

 \subsection{Special values for $\tau$\ and
$\tau$-semistability}\label{critical-values}

 In principle $\tau$\ is a continuously varying real parameter.  The stability
properties of a given triple do not likewise vary continuously, but can change
only at certain rational values of $\tau$.  This is the same phenomenon as
appears in the case of stable pairs.  In both cases it is due to the fact that,
except for $\tau$\ itself, all numerical quantities in the definition of
stability are rational numbers with bounded denominators.  In the case of
holomorphic pairs, this has the additional consequence that for the generic
choice of $\tau$\ there is no distinction between stability and semistablity.
This is in contrast to the case of pure bundles, where the notions of stability
and semistability coincide only when the rank and degree of the bundle are
coprime. The next proposition shows that for a holomorphic triple both  the
value of $\tau$\ and the greatest common divisor of the rank and degree, are
relevant.

 \begin{prop}\label{critical}Let \ $T=(E_1,E_2,\Phi)$\ be a $\tau$-semistable
triple, and let $T'=(E'_1,E'_2,\Phi')$\ be a subtriple such that
$\theta_{\tau}(T')=0.$
 Then either
 \be
 r_1r_2'=r_2r_1'\;\; \text{and}\;\; \mu(E'_1\oplus E'_2)=\mu(E_1\oplus
E_2),\label{3.15a}
 \ee
 or
 \be
\frac{r_2(r_1'+r_2')\mu(T')-r_2'(r_1+r_2)\mu(T)}{r_2r_1'-r_1r_2'}=
\tau.\label{3.
15b}
 \ee

 In particular, if \  $r_1+r_2$\ and $d_1+d_2$\ are coprime, and $\tau$\ is  a
not rational number with denominator of magnitude less than $r_1r_2$, then all
$\tau$-semistable triples are \ts.
 \end{prop}

 \pf
  From the definition of $\theta_{\tau}$, we see that $\theta_{\tau}(T')=0$\ is
equivalent to
 $$(\mu(E'_1\oplus E'_2)-\frac{r_2'}{r_2}\frac{r_1+r_2}{r_1'+r_2'}\mu(E_1\oplus
E_2))=
 \tau\frac{r_1'r_2-r_1r_2'}{r_1'r_2+r_2'r_2}.$$
 If $r_1'r_2-r_1r_2'\ne 0$\ we get (\ref{3.15b}), and if $r_1'r_2-r_1r_2'=0$\
then
 $$
 \frac{r_2'}{r_2}\frac{r_1+r_2}{r_1'+r_2'}=1
 $$
  and we get (\ref{3.15a}).
 \qed

 Next we compare the stability conditions for a triple and for the two bundles
in the triple.

 \begin{prop}\label{ts-ss}Let $(E_1,E_2,\Phi)$\ be a non-degenerate holomorphic
triple.  There is an $\epsilon>0$, which depends only on the degrees and ranks
of $E_1$\ and $E_2$, and such that for $\mu(E_1)<\tau<\mu(E_1)+\epsilon$\ the
following is true:
 \begin{tabbing}
 \ \ \ \= (0)\ \= \kill
 \>  {\em (1)}\> If $(E_1,E_2,\Phi)$\ is a $\tau$-stable triple, then both
$E_1$\ and $E_2$\ are semistable bundles.\\
 \> {\em (2)}\>  Conversely, if $E_1$\ and $E_2$\ are stable bundles, then
$(E_1,E_2,\Phi)$\ will be a $\tau$-stable\\
 \> \>  triple for any choice of $\Phi\in H^0(\Hom(E_2,E_1))$.
 \end{tabbing}
 \end{prop}
 \pf For all subbundles $E_1'\subset E_1$\ the slope $\mu(E_1')$\ is a rational
number with denominator
 less than $r_1$.  Clearly, if we pick $\epsilon$\ small enough then the
interval $(\mu(E_1),\mu(E_1)+\epsilon)$\ contains no rational numbers with
denominator less than $r_1$.  The condition $\mu(E_1')<\tau$\ is thus
equivalent to the condition $\mu(E_1')\le\mu(E_1)$, i.e. to the semistability
of $E_1$.

 Furthermore, as noted above, if $\tau < \mu(E_1)+\epsilon$\ then $\tau' >
\mu(E_2)-\frac{r_1}{r_2}
 \epsilon$. Hence if $\frac{r_1}{r_2}\epsilon$\ is small enough, then the
condition
 $\mu(E_2/E_2')>\tau'$\ for all subbundles $E_2'\subset E_2$\ becomes
equivalent to the condition that $\mu(E_2/E_2')\ge\mu(E_2)$.

 Conversely, suppose $\tau =\mu(E_1)+\delta$\ for some $\delta >0$, and that
$\Phi$\ is any section of $H^0(\Hom(E_2,E_1))$.  Then for any subtriple
$(E'_1,E'_2,\Phi)$\ we get

$$(r_1'+r_2')\theta_{\tau}\stri=r_1'(\mu(E'_1)-\mu(E_1))+r_2'(\mu(E'_2)-
\mu(E_2))+(r_2r_1'-r_1r_2')\delta,$$
 where $r_1'=\rank E'_1$\ and $r_2'=\rank E'_2$.
 If $E_1$\ and $E_2$\ are stable, and $\delta$\ is small enough, then it
follows from this that $\theta_{\tau}\stri<0$\ for all subtriples.
 \qed

\section{Main theorem}\label{theorem}

 In this section we shall show how the stability of a holomorphic triple over
$X$ relates to the stability of the associated (SU(2)-equivariant) bundle over
$\XxP$.  As in Section
 \ref{background},
 let $F\lra\XxP$\ be the extension associated to the triple $\tri$, i.e. let
$F$ be
 \be
 0\lra\ps E_1\lra F\lra \ps E_2\otimes\qs\cO(2)\lra 0,
 \label{extension-b}
 \ee
 where $p$ and $q$ are the projections from $\xp$ to $X$ and $\dP^1$
 respectively, and $\cod$ is the line
 bundle of degree 2 over $\dP^1$.
 To relate the $\tau$-stability of $\tri$ to the  stability of $ F$ we
 need to consider some \kahler\
 polarization on $\xp$. The parameter $\tau$ will be encoded in this
 polarization.
 Let us  choose a metric on $X$ with \kahler\ form $\omega_X$, with volume
normalized to one.
 The metric we shall consider on
 $\xp$ will be, as in Section \ref{background}, the product of a
 the metric on $X$ with a coefficient depending on  a parameter $\sigma>0$,
and the
 Fubini--Study metric on $\dP^1$ with volume also normalized to one. The \kah\
form
 corresponding to this metric depending on the parameter $\sig$ is
 \be
 \os=\frac{\sigma}{2}\ps\omega_X\oplus \qs\omega_{\dP^1}.\label{kah-pol}
 \ee

 We can now state the main result of this section.

 \begin{thm}\label{tsvs}
 Let $\tri$ be a  holomorphic triple over a \crs\
 $X$. Let $ F$ be the holomorphic bundle over $\xp$ defined by $\tri$ as in
Proposition \ref{equi-hvb}, and
 let
 \be
 \sig(\tau)=\frac{(r_1+r_2)\tau-(\deg E_1 +\deg  E_2)}{r_2}.\label{s-t'}
 \ee

 Suppose that in $\tri$\ the two bundles $E_1$\ and $E_2$\ are not isomorphic.
Then $\tri$ is \ts\ (equivalently $\sigma$-stable) if and only if  $ F$ is
stable
 with respect to $\os$.

 In the case that $E_1\cong E_2\cong E$, the triple $(E,E,\Phi)$ is
$\tau$-stable
 (equivalently $\sigma$-stable) if and only if  $F$ decomposes as a direct sum
 $$
 F=\ps E\otimes \qs\cO(1)\oplus \ps E\otimes \qs\cO(1),
 $$
 and $\ps E\otimes \qs\cO(1)$\ is stable with respect to $\os$.
 \end{thm}
 \pf
 As mentioned in \S\ref{stability} , if $ E_2$ is a line bundle the
 $\tau$-stability of $\tri$ is equivalent to the
 $\tau$-stability of the pair $( E_1\otimes E_2^\ast,\Phi)$ in the sense of
 Bradlow \cite{B2}. In this case
 Theorem \ref{tsvs} has been  proved
 in \cite[Theorem 4.6.]{GP3}. The main ideas of that proof extend to the
general case
 in a rather straightforward manner.

 Recall that the bundle $ F$ associated to $\tri$ comes equipped with
 a holomorphic action of $\su$. It makes sense therefore to talk about
 the \sui\ stability of $ F$. As explained in \S\ref{background}, this is
 like ordinary stability, but the slope condition has to be satisfied
 only for \sui\ subsheaves of $ F$.
 In order to prove the theorem we shall prove first the following
 slightly weaker  result.

 \begin{prop}\label{tsvis}
 Let $\tri$ be a  holomorphic triple over a \crs\
 $X$. Let $ F$ be the holomorphic bundle over $\xp$ defined by $\tri$  and
 let $\sigma$ and $\tau$ be related by ({\em \ref{s-t'}}).
 Then $\tri$ is \ts\ (equivalently $\sigma$-stable) if and only if  $ F$ is
 $\su$-invariantly stable with respect to $\os$.
 \end{prop}
 \pf
 We saw in \S\ref{background} (Lemmas \ref{equi-sub} and \ref{push-out}) that
there is a
 one-to-one correspondence
 between subtriples $T'=\stri$ of $T$ and \sui\ coherent subsheaves $F'\subset
F$.
 Moreover, the subsheaf $F'$ defined by $T'$ is an extension of the form
 \be
 \sextn.\label{subextension}
 \ee

 In terms of the parameter $\tau'$ as defined in (\ref{3.5}), the relation
 between $\sig$ and $\tau$, given by (\ref{s-t'}), can be rewritten as
 $$
 \sig=\tau-\tau'.
 $$
 If $\tri$ is \ts\ (equivalently, $\sigma$-stable) it follows from (3) in
Proposition \ref{l-bound} that $\sigma$
  defined by (\ref{s-t'}) is positive.
 The slope of $ F'$ with respect to $\os$ is defined as
 $$
 \ms( F')=\frac{\deg_\sig F'}{\rank F'},
 $$
 where $\deg_\sigma F'$ is the degree of $F'$, which is given by
 $$
 \deg_\sig F=\frac{1}{2}\int_{\xp}c_1( F)\wedge\os.
 $$
 The proposition follows now from the following lemma.
 \begin{lemma}\label{slope}Let $T'$ be a subtriple of $T$ and $F'$ the
corresponding
 \sui\ subsheaf of $F$. Let $\sig$ be as in Proposition \ref{tsvis}.
 The following are equivalent
 \begin{tabbing}
 \ \ \ \= (0)\ \= \kill
  \> {\em (1)}\> $\ms( F')<\ms( F)$\\
 \> {\em (2)}\> $\TT(T')<0$\\
 \> {\em (3)}\> $\ms(T')<\ms(T)$.
 \end{tabbing}
 \end{lemma}
 \pf
 The equivalence between (2) and (3) corresponds, of course, to the two
equivalent
 definitions of stability for $T$ (cf. Proposition \ref{tstab-sstab}).

 From (\ref{extension-b}) and (\ref{subextension}) we obtain that
 $$
 \ms( F)=\frac{\deg E_1+\deg E_2+\sigma r_2}{r_1+r_2},
 $$
 and
 $$
 \ms( F')=\frac{\deg E'_1+\deg E'_2+\sigma r_2'}{r_1'+r_2'},
 $$
 where $r_1'=\rank E_1'$ and $r_2'=\rank E_2'$.
 From Definition \ref{t-stab} we  immediately obtain the
 equivalence between (1) and (3).
 \qed

\noindent{\em Remark.}
As usual, in order for $F$ to be \suis\ it is enough to check condition (1) of
Lemma
\ref{slope} only for saturated \sui\ subsheaves, that is \sui\ subsheaves $F'$
such that
the quotient $F/F'$ is torsion-free. Such a subsheaf $F'$ is a subbundle
outside
of a set of codimension greater or equal than 2. Hence by $\su$-invariance one
concludes
that $F'$ must be  actually a subbundle of $F$ over the whole $\xp$. It is
easy to see that the saturation of $F'$ implies that of $\ps E_1'$ and
$\ps E_2'\otimes\qs\cod$ in (\ref{subextension}), and hence $E_1'\subset E_1$
and
$E_2'\subset E_2$ are in fact subbundles. In other words, the one-to-one
correspondence
between \sui\ subsheaves of $F$ and subtriples of $T$ established in Lemma
\ref{push-out}
 sends saturated subsheaves into saturated subtriples.

 To prove the theorem, we first observe that if $F$ is \suis\ then, from
Theorem \ref{isvs},
 it decomposes as a direct sum
 \be
  F= F_1\oplus F_2\oplus ...\oplus F_k \label{decomposition}
 \ee
 of stable bundles, where $ F_i$ is the transformed by an element of
 $\su$ of a fixed subbundle $ F_1$ of $ F$.

 For the remaining parts of the theorem, the proof splits into two cases,
corresponding to whether $E_1$\ and $E_2$\ are isomorphic (as holomorphic
bundles) or not.  We treat the non-isomorphic case first.  Notice that in this
case, the map $\Phi$\ certainly cannot be an isomorphism.
 Clearly if $F$ is stable it is in particular \suis\ and hence by the previous
Proposition,
 the corresponding triple will be \ts.
 Suppose now that $\tri$ is \ts, and that $\Phi$\ is not an isomorphism.
 Our strategy to prove  the stability of $ F$ will be
 to prove that $ F$ is simple, that is
 $H^0(\End F)\cong \dC$, and hence there must be just one summand
 in the decomposition of $F$ given by (\ref{decomposition}).

 To compute $H^0(\End F)\cong H^0( F\otimes F^\ast)$ let us tensor
 (\ref{extension-b}) with $ F^\ast$. We obtain the
 short exact sequence
 $$
 0\lra\ps  E_1\otimes F^\ast\lra F\otimes F^\ast\lra \ps
 E_2\otimes\qs\cod\otimes F^\ast\lra 0,
 $$
 and the corresponding  sequence in cohomology
 \be
 0\lra H^0(\ps E_1\otimes F^\ast)\lra H^0( F\otimes F^\ast)\lra
 H^0(\ps E_2\otimes\qs\cod\otimes F^\ast)\lra.
 \label{simple}
 \ee

 We first compute $H^0(\ps E_1\otimes F^\ast)$. Dualizing (\ref{extension-b}),
 tensoring with $\ps E_1$, and using that
 $H^0(\ps( E_1\otimes E_2^\ast)\otimes\qs\cO(-2))=0$, we have the
 sequence in cohomology
 \be
 0\lra H^0(\ps E_1\otimes F^\ast)\lra H^0(\ps( E_1\otimes
 E_1^\ast))\stackrel{g}{\lra}
 H^1(\ps( E_1\otimes E_2^\ast)\otimes \qs\cO(-2)).
 \label{cs1}
 \ee
 By the K\"unneth formula
 $$
 H^0(\ps( E_1\otimes E_1^\ast))\cong H^0( E_1\otimes E_1^\ast)
 \;\;\;\;\mbox{and}\;\;\;\;
 H^1(\ps (E_1\otimes E_2^\ast)\otimes\qs\cO(-2))\cong H^0( E_1\otimes
 E_2^\ast).
 $$
Now, thanks to these isomorphisms, $g$ can be interpreted as the map $H^0(
E_1\otimes E_1^\ast)\ra H^0( E_1\otimes
 E_2^\ast)$ defined by $\Phi$, i.e.
 $g(u)=u\Phi$.

 Now, from the $\tau$-stability of $\tri$ one has from Corollary
\ref{stab-simple} that
 $\tri$\  is simple.
 Thus $\Ker g\cong 0$ and from the
 exactness of (\ref{cs1}) one obtains
 \be
 H^0(\ps E_1\otimes F^\ast)=0.
 \label{van1}
 \ee
 To compute $H^0(\ps E_2\otimes \qs\cod\otimes F^\ast)$, we dualize
 (\ref{extension-b}) and tensor it with
 $\ps E_2\otimes\qs\cod$, to
 get the sequence
 \be
 0\lra H^0(\ps( E_2\otimes E_2^\ast))\lra H^0(\ps
 E_2\otimes\qs\cod\otimes F^\ast)
 \lra H^0(\ps( E_1^\ast\otimes E_2)\otimes\qs\cod).
 \label{cs2}
 \ee
 \begin{lemma}
 Let $\tri$ be \ts\ and suppose that $\Phi$ is not an isomorphism,
 then $H^0( E_1^\ast\otimes E_2)=0$.
 \label{van2}
 \end{lemma}
 \pf
 Suppose that there is a non-zero homomorphism $\Psi:  E_1\ra E_2$.
 Let $u=\Phi\Psi\in H^0(\End E_1)$
 and $v=\Psi\Phi\in H^0(\End E_2)$. Then $u\Phi=\Phi v$ and, since
 $\tri$ is simple, we have
 that

 $$
 u=\lambda \bold I_{ E_1}\;\;\;\;\mbox{and}\;\;\;\; v=\lambda \bold
I_{E_2},\;\;\;\mbox{for}
 \;\;\;\lambda\in\dC.
 $$
  If $\lambda\neq 0$, we easily see that
 $\Phi$ is an isomorphism, contradicting the assumption of the Lemma.
 Thus $\lambda=0$ and then
 $$
 \im \Psi\subseteq\Ker \Phi
 \;\;\;\;\;\mbox{and}\;\;\;\;\;
 \im \Phi\subseteq\Ker \Psi.
 $$
 We can therefore consider  the subtriples of $\tri$

 $$
 T_1=( K, E_2,\Phi)\;\;\;\;\mbox{and}\;\;\;\;T_2=(0, I,\Phi),
 $$
 where $ K=\Ker \Psi$ and $ I=\im \Psi$.
 Let $r_1'=\rank  K$, $r_1''=\rank  I$, $d_1'=\deg  K$ and $d_1''=\deg I$.
 Applying the $\tau$-stability condition to $T_1$ and $T_2$ we get the
 inequalities
 $$
 \begin{array}{l}
 r_2d_1''-r_1''(d_1+d_2)+r_1r_1''\tau<0\\
 d_1'-d_1+(r_1-r_1')\tau<0.
 \end{array}
 $$
 From this  and using that from
 $$
 0\lra K\lra E_1\lra I\lra 0,
 $$
 $r_1=r_1'+r_1''$ and $d_1=d_1'+d_1''$,
 we obtain that
 $$
 \tau<\mu( E_1\oplus E_2),
 $$
 which is equivalent to $\sig(\tau)<0$, contradicting the
 $\tau$-stability of $\tri$.
\qed

 From the K\"unneth formula and  Lemma \ref{van2} we get that
 $$
 H^0(\ps( E_1^\ast\otimes E_2)\otimes\qs\cod)\cong
 H^0( E_1^\ast\otimes E_2)\otimes H^0(\cod)\cong 0,
 $$
 and  from (\ref{cs2})

 $$
 H^0(\ps E_2\otimes\qs\cod\otimes F^\ast)\cong H^0( E_2\otimes
 E_2^\ast).
 $$
 From this  and (\ref{van1}) the first three terms in  (\ref{simple})
 reduce to
 $$
 0\lra H^0( F\otimes F^\ast)\stackrel{i}{\lra} H^0( E_2\otimes
 E_2^\ast).
 $$
 Since $ F$ is \suis\ then it is $\su$-invariantly simple, i.e. the only \sui\
 endomorphisms are  multiples of the identity.
 Let $\Psi\in H^0( F\otimes F^\ast)$ be a non \sui\ endomorphism of $
 F$, i.e. $\Psi^g\neq\Psi$ for some
 $g\in \su$. Since $i$ must be compatible with the action of $\su$, we
 get
 $$
 i(\Psi^g)=(i(\Psi))^g.
 $$
 On the other hand, since the action of $\su$ on $H^0( E_2\otimes
 E_2^\ast)$ is trivial
 $$
 (i(\Psi))^g=i(\Psi),
 $$
 hence $ i(\Psi)=i(\Psi^g)$ contradicting the injectivity of $i$. Thus
 $H^0( F\otimes F^\ast)\cong\dC^\ast$,
 which concludes the proof of our theorem for the case where $E_1$\ and $E_2$\
are not isomorphic.

 Now suppose that $E_1\cong E_2$.  We first prove

 \begin{lemma}\label{ts-iso} Suppose $E_1\cong E_2$.  Then for any
 $\tau>\mu(E_1)$, the triple $(E_1, E_2,\Phi)$\ is $\tau$-
 stable if and only if $\Phi$\ is an isomorphism and $E_1$\
 is  stable.
 \end{lemma}
 \pf
 Suppose that the triple is $\tau$-stable.  Then (by Corollary 3.8)
 $\Phi$\ is an isomorphism.  Now consider the subtriples of
 the
 form $T'=(\Phi(E'_2), E'_2,\Phi)$.  These have $r_1'=r_2'$\ and,
 since $\Phi$\ is an isomorphism, $\mu(T')=\mu(E'_2)$.
 Hence $\theta_{\tau}(T')=\mu(E'_2)-\mu(E_2)$, and thus the
 $\tau$-stability of the triple implies the stability of
 $E_2$.
 Conversely, if $E_2\cong E_1$\ and both are stable, then
 all non trivial $\Phi$ are isomorphisms.  It now follows as a
 special case of Proposition 3.6 that the triple $(E_1, E_2,\Phi)$\
 is $\tau$-stable.
 \qed

 Suppose that $E_1\cong E_2\cong E$\ and  that the bundle $F$ associated to
 $(E,E,\Phi)$ is of the form
 \be
 F=\ps E\otimes\qs\cO(1)\oplus \ps E\otimes\qs\cO(1).\label{2-sum}
 \ee
 If we assume now that  $\ps E\otimes \qs\cO(1)$ is  stable, then $E$ is also
stable and
 hence
 $ H^0(E\otimes E^*)\cong\dC$. Thus  there is only one non-trivial
 extension class (corresponding to
 $\Phi=\lambda \bold I$).  We must now examine this (unique) non-trivial
extension
 $$
 0\lra \ps E \lra F \lra \ps E\otimes q^*\cal O(2)\lra 0.
 $$
 This is of course nothing else but the pull-back to $\xp$ of the non-trivial
 extension on $\dP^1$
 $$
 0\lra\cO\lra\cO(1)\oplus\cO(1)\lra\cO(2)\lra 0,
 $$
 tensored with $\ps E$. Thus
 the action of $\su$ permutes the two summands in (\ref{2-sum}) and from
Theorem \ref{isvs}
 we conclude that $F$ is an \suis\ bundle. The $\tau$-stability of $\tri$
follows
 now from Proposition \ref{tsvis}.

 Conversely, suppose that $\tri$ is \ts, then from Lemma \ref{ts-iso} we obtain
that
 $E_1\cong E_2\cong E$ is stable and $\Phi$ is hence a non-trivial multiple of
the identity.
 From the above discussion we conclude that
 $$
 F=\ps E\otimes\qs\cO(1)\oplus \ps E\otimes\qs\cO(1).
 $$
 On the other hand, from Proposition \ref{tsvis}, we
  argue as before that $F$ \ is certainly invariantly
 stable . Also (\ref{simple}), (\ref{cs1}) and (\ref{van1}) show that we have
an exact
 sequence
 $$0\lra H^0(F\otimes F^*)\lra H^0(p^*E\otimes q^*\cal O(2)\otimes F^\ast).$$
 Using that  $H^0(E\otimes E^*)\cong \Bbb C$, the exact sequence (\ref{cs2})
becomes,
 $$0\lra \Bbb C\lra H^0(p^*E\otimes q^*\cal O(2)\otimes F^*)
 \lra \Bbb C^3.$$

 From these two exact sequences, we see that
 \be
 1\le h^0(F\otimes F^*)\le h^0(p^*E\otimes q^*\cal O(2)\otimes F^*)\le 4\ .
 \label{constraint}
 \ee

 But since $F$\ is invariantly stable it is given by the direct sum
 (\ref{decomposition}).
   In that case, $$H^0(F\otimes F^*)\cong GL(k,\Bbb C)\ ,$$
 where $k$\ is the number of stable summands in $F$.  It follows that
 since  $h^0(F\otimes F^*)\ne 1$, then $h^0(F\otimes F^*)=k^2-1$\ for some
integer $k$.
 The only possibility consistent with the constraint (\ref{constraint}) is thus
  $h^0(F\otimes F^*)=3$, i.e. $k=2$. Hence the bundle $\ps E\otimes\qs\cO(1)$
in the decomposition of $F$ is stable, which finishes the proof of the Theorem.
 \qed

 Notice  that the conclusion of Proposition \ref{tsvis}  extends
straightforwardly to cover
 polystable objects.  We thus get the following corollary, which will be useful
in the next
 section

 \begin{thm}\label{ptsvps}
 Let $\tri$ be a holomorphic triple and $F$ be the corresponding holomorphic
 bundle over $\xp$.  Let $\tau$\ and $\sigma$\ be related as above. Then $\tri$
is $\tau$-polystable if and only if
 $F$ is $\su$-invariantly polystable with respect to $\os$.
 \end{thm}
 \section{Relation to vortex equations}
 In this section we relate the $\tau$-stability of a triple
 to the existence of solutions to the coupled $\tau$-vortex
 equations.  Using the idea of dimensional reduction this can
 be viewed as simply a special case of the Hitchin-Kobayashi
 correspondence between stable bundles and Hermitian-Einstein
 metrics.  Indeed, specializing to the case of SU(2)-equivariant bundles over \
 $\xp$\ as in \S \ref{background} (i.e. to the case of the bundles associated
with triples over $X$), we have already  assembled seen the results we need.

 From \cite{GP3} (cf. also Theorem \ref{inv-h-k}) we know that the
 SU(2)-equivariant bundle $F$\ admits a Hermitian-Einstein metric
 with respect to $\omega_{\sigma}$\ if and only if
 $F$\ is polystable with respect to $\omega_{\sigma}$.
 By Theorem \ref{ptsvps} we know that $F$\ is $\su$ invariantly polystable with
 respect to $\omega_{\sigma}$\ if and only if the corresponding triple
 $\tri$\ is $\tau$-polystable with $\tau$\ and $\sigma$\ related by
(\ref{s-t'}).  Finally,
 by Proposition \ref{dr} there is a bijective correspondence
 between the SU(2)-equivariant Hermitian-Einstein
 metrics on $F$ and the solutions to the coupled vortex
 equations (\ref{cves}) on $\tri$.

 Putting all this together, we see that we have filled in
 three sides of the following ``commutative diagram'' for a
 holomorphic triple over $X$\ and the corresponding
 SU(2)-equivariant bundle over $\xp$\ (with \kahler\ form $\omega_{\sigma}$)

\begin{tabular}{cccc}
 & & & \\
 & SPECIAL &\ \  Hitchin-Kobayashi \ \ & HOLOMORPHIC \\
 & METRICS & correspondence & INTERPRETATION \\
 & & & \\
 &   SU(2)-invariant Solutions & Th. 2.18  & SU(2)-invariantly \\
On $\xp$ & to Hermitian-Einstein & & Polystable \\
 & Equations & & Extensions \\
 & & & \\
dim.\ red. &Prop. 2.14  & &Th. 4.6 \\
& & &\\
 & Solutions to & & $\tau$-stable\\
On $X$ & Coupled $\tau$-Vortex & & holomorphic\\
 & Equations & & Triples\\

\end{tabular}

 By tracing around the three sides of this diagram which
 already have arrows filled in, we can fill in the arrows on
 the fourth side and prove

 \begin{thm}\label{existence}  Let $T=(E_1,E_2,\Phi)$\ be a holomorphic triple.
 Then the
 following are equivalent.

 \ \ \ {\em (1)}\ The bundles support Hermitian metrics $h_1, h_2$\
 such that the coupled $\tau$-vortex equations
 are satisfied, i.e. such that
 \be
  \sqrt{-1}\Lambda F_{h_1}+\Phi \Phi^*=2\pi\tau\bold I_{E_1},\label{4.1a}
 \ee
 \be
 \sqrt{-1}\Lambda F_{h_2}-\Phi^*\Phi=2\pi\tau'\bold I_{E_2},\label{4.1b}
 \ee
 with
 \be
 r_1\tau+r_2\tau'=\deg E_1 + \deg E_2,
 \label{4.1c}
 \ee

 \ \ \ {\em (2)}\ The triple is $\tau$-polystable.
 \end{thm}

 Notice that the statement and conclusion of this theorem
 make no mention
 of $\xp$\ or the \sue\ bundle $F$.  One might
 thus expect a
 more direct proof that does not use dimensional reduction.
 We will not
 attempt to prove both directions of the biconditional in the
 theorem in
 this way, but the one direction is quite
 easily seen.
 That is, one can show how the $\tau$-stability condition
 can be derived directly as a consequence of the coupled
 vortex  equations.   We do this as follows.

 Let $T'=(E'_1,E'_2,\Phi)$\ be a holomorphic saturated sub-triple of
 $T$, and let
 $$E_1=E'_1\oplus E''_1,$$
 and
 $$E_2=E'_2\oplus E''_2$$
 be smooth orthogonal splittings of $E_1$\ and $E_2$.  With
 respect to these splittings, we get a block diagonal
 decomposition of $\sqrt{-1}\Lambda F_{h_1}$, $\sqrt{-1}\Lambda
 F_{h_2}$as

 \be
 \sqrt{-1}\Lambda F_{h_i}=\left(\begin{array}{cc}
 \sqrt{-1}\Lambda F'_i+\Pi_i & \ast\\
 \ast& \sqrt{-1}\Lambda F''_i-\Pi_i
 \end{array}\right),\label{4.2}
 \ee
 where $\Pi_i$\ is a positive definite endomorphism coming
 from the second fundamental form for the inclusion of
 $E'_i$\ in $E_i$.  We also get a decomposition of $\Phi$\ as
 \be
 \Phi=\left(\begin{array}{cc}
 \Phi' & \Theta\\
 0 & \Phi''
 \end{array}\right).\label{4.3}
 \ee
 The coupled vortex equations thus split into equations on
 the summands of $E_1$\ and $E_2$\ to yield
 \begin{eqnarray}
  \sqrt{-1}\Lambda F'_1+\Pi_1+\Phi' (\Phi')^*+\Theta\Theta^*
 &=2\pi\tau\bold I'_1,\label{4.4a}\\
 \sqrt{-1}\Lambda F''_1-\Pi_1+\Phi''
 (\Phi'')^*&=2\pi\tau\bold I''_1,\label{4.4b}\\
 \sqrt{-1}\Lambda F'_2+\Pi_2-(\Phi')^*
 \Phi'&=2\pi\tau'\bold I'_2,\label{4.4c}\\
 \sqrt{-1}\Lambda F''_2-\Pi_2-(\Phi'')^* \Phi''-
 \Theta^*\Theta &=2\pi\tau'\bold I''_2.\label{4.4d}
 \end{eqnarray}

 We now integrate the trace of these equations over the base
 manifold $X$.  We use the notation
 $$
 d'_i=\deg E'_i=\frac{\sqrt{-1}}{2\pi}\int_X \tr(\Lambda F'_i),
 $$
 $$
 d''_i=\deg E''_i=\frac{\sqrt{-1}}{2\pi}\int_X \tr(\Lambda F''_i),
 $$
 $$
 ||\Pi_i||^2=\frac{1}{2\pi}\int_X \tr(\Pi_i),
 $$
 $$
 ||\Theta||^2=\frac{1}{2\pi}\int_X \tr(\Theta\Theta^*),
 $$

 From equations (\ref{4.1a}) and (\ref{4.1c}) we get
 \be
 d''_1+d''_2 -
(||\Pi_1||^2+||\Pi_2||^2)-||\Theta||^2=\tau(r_1-r_1')+\tau'(r_2-r_2'),
 \label{4.5}
 \ee
 and from (\ref{4.1b}) and (\ref{4.4a}) we get
 \be
 d'_1+d'_2 + (||\Pi_1||^2+||\Pi_2||^2)+||\Theta||^2=\tau r_1'+\tau'
r_2'.\label{4.6}
 \ee
 Using the fact that $d_i=d'_i+d''_i$, these equations can be
 combined to give
 \be
 \deg(E'_1\oplus
E'_2)+2(||\Pi_1||^2+||\Pi_2||^2+||\Theta||^2)=r_1'\tau+r_2'\tau'\label{4.7}
 \ee
 It now follows by the positivity of the terms $||\Pi_i||^2$\
 and $||\Theta||^2$, that
 \be
\mu(T')\le \frac{r_1'\tau+r_2'\tau'}{r_1'+r_2'}.
\label{4.8}
\ee
 But from equation (\ref{4.1c}) we have that
 $$\frac{r_1\tau+r_2\tau'}{r_1+r_2} = \mu(T),$$
 and by using this to solve for $\tau'$, one sees that
 (\ref{4.8}) is equivalent to
 the condition $\theta_{\tau}(T')\le 0$.  Furthermore, it
 follows from (\ref{4.7}) that in order to have $\theta_{\tau}(T')=
 0$, one needs
 $$||\Pi_1||^2=||\Pi_2||^2=||\Theta||^2=0.$$
 This means that the bundles split holomorphically as
 $E_i=E'_i\oplus E''_i$\ and $\Phi=\Phi'\oplus \Phi''$, i.e.
 the triple $T=\tri$\ splits as $T=(E'_1,E'_2,\Phi')\oplus
 (E''_1,E''_2,\Phi'')$.  Each summand separately supports a
 solution to the coupled $\tau$-vortex equations.  It is
 possible that in this splitting $\Phi'$\ or $\Phi''$\ is
 trivial.  The bundles in the degenerate subtriple then each
 supports solutions to the Hermitian-Einstein equations.  We
 have thus proven the following proposition:

 \begin{prop}Let $T=(E_1,E_2,\Phi)$\ be a
 holomorphic triple in which the bundles support Hermitian
 metrics $h_1, h_2$\ such that the coupled $\tau$-vortex
 equations are satisfied.  Then
 \begin{tabbing}
 \ \ \ \=(0)\ \= \kill
 \> {\em (1)}\> $T$\ splits as a direct sum of triples
 $(E_{1i},E_{2i},\Phi_i)$,
 i.e. $E_1=\bigoplus E_{1i}$, $E_2=\bigoplus E_{2i}$,\\
 \> \> and $\Phi=\bigoplus \Phi_i$,\\
 \> {\em (2)}\> each summand $(E_{1i},E_{2i},\Phi_i)$\ is
 either a $\tau$-stable triple, or $\Phi_i=0$\ and both \\
 \> \> bundles are stable. In the latter case, the slope of $E_{1i}$\
 is $\tau$\ and the slope of $E_{2i}$\\
\> \> is $\tau'$.
 \end{tabbing}
 That is $T$\ is $\tau$-polystable.
 \end{prop}
\section{Moduli spaces}\label{moduli}
\subsection{Moduli spaces of stable triples}
Recall that two triples $T=\tri$ and $T'=\stri$ are isomorphic if there exist
isomorphisms
$u:E_1\lra E_1'$\ and $v:E_2\lra E_2'$\ making the following
diagram commutative
$$
\begin{array}{ccc}
\; E_2& \stackrel{\Phi}{\lra} & \; E_1 \\
 v{\downarrow}&  &  u{\downarrow}\\
\; E_2'& \stackrel{\Phi'}{\lra} & \;  E_1'.
\end{array}
$$

After fixing the topological invariants of our bundles, that is the ranks $r_1$
and $r_2$
and the first Chern classes $d_1$ and $d_2$, let $\glM$ be the set of
equivalence classes
of holomorphic triples and $\glM_\tau\subset\glM$ be the subset of equivalence
classes of
\ts\ triples. Our goal in this section is to show that $\glM_\tau$ has the
structure of
an algebraic variety, more precisely:

\begin{thm}\label{MODULI}
Let $X$ be a \crs\ of genus $g$ and let us fix ranks $r_1$ and $r_2$ and
degrees $d_1$
and $d_2$. The moduli space of \ts\ triples $\glM_\tau$ is a complex analytic
space
with a natural \kahler\ structure outside of the singularities. Its dimension
at a smooth
point is
\be
1+r_2d_1-r_1d_2+(r_1^2+r_2^2-r_1r_2)(g-1).
\label{dimension}
\ee
The moduli space, $\glM_\tau$ is non-empty if and only if $\tau$ is inside the
interval
\be
(\mu(E_1), \mu_{MAX})
\label{interval}
\ee
where

\be
\mu_{MAX}=\mu(E_1)+\frac{r_2}{|r_1-r_2|}(\mu(E_1)-\mu(E_2))
\ee
if $r_1\neq r_2$,  and $\mu_{MAX}=\infty$\ if $r_1=r_2$.

Moreover $\glM_\tau$ is in general a quasi-projective variety.  It is in fact
projective if
$r_1+r_2$ and $d_1+ d_2$ are coprime and $\tau$ is generic.
\end{thm}
\pf
There are several approaches one can  take to prove this theorem. One can use
standard Kuranishi deformation methods as done in \cite{B-D1,B-D2} for the
construction of the moduli spaces of stable pairs. Alternatively one can use
geometric
invariant theory methods to give an algebraic geometric construction of our
moduli spaces,
generalizing the construction of the moduli space of stable pairs given in
\cite{Be,T}.
We will leave these two direct methods for a future occasion and instead will
exploit the relation between \ts\ triples and equivariant bundles over $\xp$.
This  method,
which  is
used in \cite{GP3} to construct the moduli spaces of triples when $E_2$ is a
line bundle,
leads also to an alternative construction of the moduli spaces of stable pairs.
Apart from smoothness considerations, which we shall discuss later, the
arguments of the proof
are the same that those in \cite{GP3}.

Let $\sigma$ be related to $\tau$ by (\ref{s-t'}) and let $\os$ be the \kahler\
form  on $\xp$ defined
by (\ref{kah-pol}). Let $\cM_\sigma$ be the moduli space of stable bundles with
respect to
$\os$ whose
underlying smooth bundle is defined by (\ref{equi-f}).
Let us exclude for the moment the case $r_1=r_2$
and $d_1=d_2$. Let $F\lra\xp$ be the bundle associated to $\tri$ as in
Proposition
\ref{equi-hvb}.
 Theorem \ref{tsvs} says that the correspondence $\tri\longmapsto F$ defines a
map
$$
\glM_\tau\lra\cM_\sigma.
$$
The action of $\su$ on $\xp$ defined in Section  \ref{background} induces an
action on
$\cM_\sigma$ and, since the bundle $F$ associated to $\tri$ is \sue\ the image
of the above map is contained in
$\cM_\sigma^{\su}$---the set of fixed points of $\cM_\sigma$ under the $\su$
action. As proved
in \cite[Proposition 5.3]{GP3} the set $\cM_\sigma^{\su}$ can be described as a
disjoint union of a finite
number of sets
$$
\cM_\sigma^{\su}=\bigcup_{i\in I}\cM_\sigma^i.
$$
The index $I$ ranges over the set of equivalence classes of different smooth
\sue\
structures on the smooth bundle $F$ defined by
(\ref{equi-f}).  Of course the way of writing $F$ in (\ref{equi-f}) already
exhibits  a
particular \sue\ structure, but in principle the bundle $F$ might
admit different ones.  The set $\cM_\sigma^i$ corresponds to the set of
equivalence
classes in $\cM_\sigma$ admitting a representative which is \sue\ for the
smooth equivariant
structure defined by $i\in I$. An equivariant smooth structure defines an
action on the space
of smooth automorphisms of the bundle $F$ and, as shown in \cite[Theorem
5.6]{GP3} the sets $\cM_\sigma^i$
can be described as the set of equivalence classes of \sue\ holomorphic
structures on
the underlying smooth \sue\  bundle defined by $i$, modulo \sue\ isomorphisms.

Let $i_0$ be the $C^\infty$ \sue\ structure on $F$ defined by (\ref{equi-f}).
As shown in Proposition
\ref{equi-hvb} there is a one-to-one correspondence
\be
\{\mbox{holomorphic triples}\}\stackrel{1-1}{\longleftrightarrow}
\{\mbox{$i_0$-equivariant holomorphic vector bundles}\}.
\label{tri-equiv}
\ee
On the other hand by Lemma \ref{equi-hom} the equivariant homomorphisms between
two equivariant holomorphic
bundles $F$ and $F'$ corresponding to triples $T$ and $T'$, respectively, are
in one-to-one
correspondence with the morphisms between $T$ and $T'$.
In fact the correspondence (\ref{tri-equiv}) descends to the quotient and thus
from Theorem \ref{tsvs} we can identify $\glM_\tau$ with $\cM_\sigma^{i_0}$.
The properties of
$\glM_\tau$ follow now from standard facts about the more familiar moduli
spaces of stable
bundles $\cM_\sigma$ \cite{D-K,G,M,Ko},
 and more particularly of the fixed-point sets $\cM_\sigma^i$
(see \cite[Theorem 5.6]{GP3} for details). Namely,

\begin{thm}\label{fix-points}$\cM_\sigma^i$ is a complex analytic variety. A
point $[F]\in \cM_\sigma^i$ is
non-singular if it is non-singular as a point of $\cM_\sigma$. The tangent
space at such a
point can be identified with the \sui\ part of $H^1(\xp, \End F)$.
$\cM_\sigma^i$ has a natural \kahler\ structure induced from that of
$\cM_\sigma$.
Moreover if $\sigma$ is a rational number then $\cM_\sigma^i$ is a
quasi-projective variety.
\end{thm}

{}From this theorem and  the identification of $\glM_\tau$ with
$\cM_\sigma^{i_0}$ we deduce that
$\glM_\tau$ is a complex analytic variety with a \kahler\ metric outside the
singularities.
To compute the dimension of the tangent space at a smooth point $[T]$ it
suffices to compute
the dimension of the \sui\ part of $H^1(\xp,\End F)$. This can be done in a
similar way  to
that of \cite[Theorem 5.13]{GP3} to obtain that
$$
\mbox{dim}\glM_\tau=1+\chi(E_1\otimes E_2^\ast)-\chi(\End E_1)-\chi(\End E_2),
$$
which by Riemann-Roch yields (\ref{dimension}).

We consider now the case $r_1=r_2=r$ and $d_1=d_2=d$. In this case by Lemma
\ref{ts-iso} we can identify
the moduli space $\glM_\tau$ with  the moduli space of stable bundles of rank
$r$ and
degree $d$ on $X$. The theorem follows now from well-known results about this
moduli space \cite{A-B,N-S}.

The fact that $\glM_\tau$ is empty outside the interval (\ref{interval}) if
$r_1\neq r_2$ and
outside $(\mu(E_1),\infty)$ if $r_1=r_2$ follows from Proposition
\ref{u-bound}.
As explained in Proposition \ref{critical} the non-generic values divide this
intervals in subintervals in such a way
that the stability properties of a given triple do not change for two values of
$\tau$ in the
same subinterval. Therefore  we can always choose $\tau$ (and hence $\sigma$)
to be
rational, which by Theorem \ref{fix-points} gives that $\glM_\tau$ is
quasi-projective.

To show the compactness of $\glM_\tau$ when $r_1+r_2$ and $d_1+d_2$ are coprime
and $\tau$ is
generic (we are also assuming that $r_1\neq r_2$ or $d_1\neq d_2$) we consider
a sequence of points in $\cM_\sigma^{i_0}$. This sequence must converge
in $\overline{\cM_\sigma}$---the Uhlenbeck compactification of $\cM_\sigma$.
Using
$\su$-invariance one can see that the limit has to correspond to a polystable
element, but
by Proposition \ref{critical} this has to be actually stable, that is the limit
must be in $\cM_\sigma$
and hence in $\cM_\sigma^{i_0}$ since this is closed. The compactness when
$r_1=r_2$ and
$d_1=d_2$ follows from the compactness of the moduli space of stable bundles of
rank $r$ and
degree $d$ when $r$ and $d$ are coprime.
The compactness of $\glM_\tau$ can  also be obtained (as it is done for pairs
in
\cite{B-D1}) from the fact that it  can
be identified with the moduli space of solutions to the coupled vortex
equations and these
 are moment map equations as we shall explain later.
\qed

It was shown in \cite[Theorem 5.13]{GP3} that when $E_2$ is a line bundle our
moduli spaces are smooth
for every value of $\tau$. This does not seem to be the case when $E_2$ is of
arbitrary
rank. However we can show the following

\begin{prop}Let $T=\tri$ be a holomorphic triple such that $\Phi$ is either
injective or surjective, then $[T]$ is a smooth point of $\glM_\tau$.
\end{prop}
\pf
Let
\be
\extn, \label{extension-3}
\ee
be the extension over $\xp$ corresponding to $T$. To prove the smoothness of
$\glM_\tau$
at the point $[\tri]$ it suffices to show that $H^2(\xp,\End F)=0$.
Tensoring (\ref{extension-3}) with $F^\ast$  the last terms in the
corresponding long exact sequence are
\be
H^2(\ps E_1\otimes F^\ast)\lra
H^2(F\otimes F^\ast)\lra H^2(\ps E_2\otimes\qs\cod\otimes F^\ast)\lra 0.
\label{2-coh-seq}
\ee
By Serre duality
$$
H^2(\ps E_1\otimes F^\ast)\cong H^0(\ps( E_1^\ast\otimes K)\otimes F)^\ast
$$
$$
H^2(\ps E_2\otimes\qs\cO(2)\otimes F^\ast)\cong
H^0(\ps (E_2^\ast\otimes K)\otimes\qs\cO(-4)\otimes F)^\ast,
$$
where $K$ is the canonical line bundle of $X$.

It is easy to see that $H^0(\ps (E_2^\ast\otimes K)\otimes\qs\cO(-4)\otimes
F)\cong 0$.
To analyse $H^0(\ps( E_1^\ast\otimes K)\otimes F)$ we tensor
(\ref{extension-3}) with
$\ps(E_1^\ast\otimes K)\otimes\qs\cO(-2)$, and since
$H^0(\ps(E_1\otimes E_1\otimes K)\otimes\qs\cO(-2))\cong 0$, we obtain
\be
0\lra H^0(\ps(E_1^\ast\otimes K)\qs\cO(-2)\otimes F)\lra
H^0(\ps(E_1^\ast\otimes E_2\otimes K))\stackrel{f}{\lra}
H^1(\ps(E_1\otimes E_1^\ast\otimes K \otimes\cO(-2)).\label{coh-sequence}
\ee
The map $f$ in the above sequence is essentially the map
$$
\begin{array}{ccc}
H^0(E_1^\ast\otimes E_2\otimes K)& \stackrel{f}{\lra}& H^0(E_1\otimes
E_1^\ast\otimes K)\\
\Psi&\longmapsto&\Phi\circ\Psi.
\end{array}
$$
Assume now that $\Phi$ is injective, if we prove that $f$ is also injective, by
the exactness
of (\ref{coh-sequence}) we would be done. Suppose that $\Ker f\neq 0$. This
means that
there exists a non-zero map $\Psi: E_1\lra E_2\otimes K$, and since
$\Psi\circ\Phi=0$,
 $\im \Psi$ is a non-trivial subsheaf contained in $\Ker \Phi$ contradicting
the injectivity.

To prove smoothness when $\Phi$ is surjective, we consider the dual triple
$T^\ast=\dtri$.
$\Phi^\ast$ is now injective and the result follows from the fact that $T=\tri$
is a smooth point if and only if $T^\ast=\dtri$ is a smooth point.

\subsection{Abel--Jacobi maps}
As shown in Proposition \ref{ts-ss} there is a range for the parameter $\tau$
such that the
$\tau$-stability of a triple $\tri$ implies the semistability of $E_1$ and
$E_2$.
Let $\glM_0$ be the moduli space of \ts\ triples for $\tau$ in such a range.
Let $N(r,d)$ be the Seshadri compactification of the moduli space of stable
bundles
of rank $r$ and degree $d$ over $X$, that is, the space of $S$-equivalence
classes
of semistable bundles.

There are natural `` Abel--Jacobi'' maps $\pi_1$ and $\pi_2$
$$
\begin{array}{ccc}
  \glM_0  &  \stackrel{ \large \pi_2}{\lra}& N(r_2,d_2)\\
          &                        &            \\
  \pi_1{\downarrow} &                &           \\
          &                        &            \\
  N(r_1,d_1)  &                    &
\end{array}
$$
defined as
$$
\pi_1([\tri])=[E_1]\;\;\;\mbox{and}\;\;\;\pi_2([\tri])=[E_2].
$$
We know also  from Proposition \ref{ts-ss} that if both $E_1$ and  $E_2$ are
stable
then the intersection of the  fibres
$\pi_1^{-1}([E_1])$ and  $\pi_2^{-1}([E_2])$ can be identified with
$\dP(H^0(E_1\otimes E_2^\ast))$. In general, though, this intersection
for non-stable points is  hard to describe.

If $\mu(E_1\otimes E_2^\ast)> 2g-2$, that is if
$$
r_2d_1-r_1d_2>r_1r_2(2g-2),
$$
where $g$ is the genus of $X$, then $H^1(E_1\otimes E_2^\ast)=0$ for $E_1$ and
$E_2$
stable and the projection from $\glM_0$ to $N(r_1,d_1)\times N(r_2,d_2)$ is a
fibration
on the the stable part.

Recall that if  $(r_1,d_1)=1$ and $(r_2,d_2)=1$ then stability and
semistability
coincide and there exist universal bundles
$$
\dE_1\lra X\times N(r_1,d_1)\;\;\;\mbox{and}\;\;\;\dE_2\lra X\times N(r_2,d_2).
$$
Let us denote by $p_1$, $p_2$ and $\pi$ the projections from
$X\times N(r_1,d_1)\times N(r_2,d_2)$ to $X\times N(r_1,d_1)$, $X\times
N(r_2,d_2)$,
and $N(r_1,d_1)\times N(r_2,d_2)$ respectively. It is clear that $\glM_0$ can
be identified with
\be
\dP(\pi_\ast(p_1^\ast\dE_1\otimes p_2^\ast \dE_2^\ast)).\label{proj}
\ee
But in the non-coprime situation we have no universal bundles $\dE_1$ and
$\dE_2$ available and the analogue of (\ref{proj}) has to be constructed as  a
moduli space in its own
right.

As explained in Theorem \ref{MODULI} the moduli space of \ts\ triples is
non-empty if and only
$\tau$ is in the interval $I=(\mu(E_1),\mu_{MAX})$. We saw in \S
\ref{critical-values} that the stability
properties of a given triple can change only at certain rational values  of
$\tau$
(the critical values) which  divide $I$ in a finite number of subintervals.
The moduli spaces for values of $\tau$ in the same open subinterval are then
isomorphic,
and they
might change only when crossing one of the critical values. We expect that, as
in the
case of stable pairs \cite{B-D-W, T}, the moduli spaces for consecutive
intervals must be
related by some sort of flip-type birational transformation. This, as well as
the construction of
a ``master'' space for triples (cf. \cite{B-D-W}) containing the moduli space
of triples
for all possible values of $\tau$, will be dealt with in a future paper.

\subsection{Vortices}
Thanks to our existence theorem the moduli space of stable triples can be
interpreted as the
moduli space of solutions to the coupled vortex equations. To understand the
meaning of this
statement one needs to regard the vortex equations as equations for unitary
connections instead of
equations for metrics. This point of view corresponds to the fact that fixing a
holomorphic
structure and varying the metric on a vector bundle is equivalent to fixing the
metric and varying the
holomorphic structure---or the corresponding connection. Recall that the space
of unitary
connections on a smooth Hermitian vector bundle can be identified with the
space of
$\dbar$-operators which in turn corresponds with the space of holomorphic
structures on our
bundle.

Let $E_1$ and $E_2$ be smooth vector bundles over $X$ and $h_1$ and $h_2$ be
Hermitian
metrics on $E_1$ and $E_2$ respectively. Let $\cA_1$ (resp. $\cA_2$) be the
space of unitary
connections on $(E_1,h_1)$ (resp. $(E_2,h_2)$).
Let $(A_1,A_2,\Phi)\in \cA_1 \times \cA_2\times \Omega^0(\Hom(E_2,E_1))$. The
vortex
equations can be regarded as the equations for $(A_1, A_2,\Phi)$
  \be
  \left. \begin{array}{l}
\dbar_{A_1\ast A_2}\Phi=0\\
 \sqrt{-1}  \Lambda F_{A_1}+\Phi\Phi^\ast=2\pi\tau \bold I_{E_1}\\
 \sqrt{-1} \Lambda F_{A_2}-\Phi^\ast\Phi=2\pi\tau'\bold I_{E_2}
  \end{array}\right \}.\label{c-cves}
  \ee
The connections $A_1$ and $A_2$ induce holomorphic structures on $E_1$ and
$E_2$
and the first equation in (\ref{c-cves}) simply says that $\Phi$ must be
holomorphic.

Let $\cG_1$ and $\cG_2$ be the gauge groups of unitary transformations of
$(E_1,h_1)$ and
$(E_2,h_2)$ respectively. $\cG_1\times\cG_2$ acts on
$\cA_1\times\cA_2\times \Omega^0(\Hom(E_2,E_1))$ by the rule
$$
(g_1,g_2).(A_1, A_2, \Phi)=(g_1A_1g_1^{-1},g_2A_2g_2^{-1},g_1\Phi g_2^{-1}).
$$
The action of \ $\cG_1\times\cG_2$\  preserves the equations and the moduli
space of
{\em coupled $\tau$-vortices} is defined as the space of all solutions to
(\ref{c-cves})
modulo this action.

 The moduli space of vortices can be obtained as a symplectic  reduction
(see \cite[Section 2.2]{GP3}) in a similar way to the moduli space of \he\
connections:
$\cA_1\times\cA_2\times\Omega^0(\Hom(E_2,E_1))$ admits a \kahler\ structure
which is
preserved by the action of $\cG_1\times\cG_2$. Associated to this action there
is a moment
map given precisely by
\be
(A_1,A_2,\Phi)\longmapsto (\Lambda
F_{A_1}-\sqrt{-1}\Phi\Phi^\ast+2\sqrt{-1}\pi\tau,
\Lambda F_{A_2}+\sqrt{-1}\Phi^\ast\Phi+2\sqrt{-1}\pi\tau').\label{mp}
\ee
Let $\mu$ be this moment map restricted to the subvariety
$$
\cN=\{(A_1,A_2,\Phi)\in \cA_1\times\cA_2\times\Omega^0(\Hom(E_2,E_1)) \;\;
|\;\;
\dbar_{A_1\ast A_2}\Phi=0\}.
$$
The moduli space of $\tau$-vortices is then  nothing else but the symplectic
quotient
$$
\mu^{-1}(0,0)/\cG_1\times\cG_2,
$$
and Theorem \ref{existence} can be reformulated by saying that there is a
one-to-one correspondence
$$
\mu^{-1}(0,0)/\cG_1\times\cG_2\stackrel{1-1}{\longleftrightarrow}\glM_\tau.
$$

\section{Some generalizations}

1. Although for simplicity we have worked  on a Riemann surface, most of our
results extend
in  a straightforward manner
to a compact complex manifold of arbitrary dimension. Of course, as in ordinary
stability, one
 needs to choose a \kahler\ metric in order to define the degree of a coherent
sheaf and
hence the slopes involved in the definition of $\tau$-stability for a triple.

2. One of our main goals in this paper  has been to show that our stability
condition for a
triple corresponds to the existence of solutions to the coupled vortex
equations. This is
the main reason for defining  our stability criterium only for vector bundles.
One can  more generally define $\tau$-stability for a triple consisting
of two (torsion free) coherent sheaves and a morphism between them. The main
results of
Sections \ref{stability} and \ref{theorem} go through in this more general
situation.

3. A. King \cite{K} has been able to characterize all \sue\ holomorphic vector
bundles
on $\xp$. Generalizing the  results in Section \ref{background}, he has shown
that these
bundles are in one-to-one correspondence with  $(2n-1)$-tuples consisting
of $n$ holomorphic vector bundles $E_1$, ..., $E_n$ over $X$ and a chain of
morphisms
$$
E_n\stackrel{\Phi_{n-1}}{\lra} E_{n-1}\lra ...\lra E_2\stackrel{\Phi_1}{\lra}
E_1.
$$
He has defined a stability condition for such a  $(2n-1)$-tuple which  involves
$(n-1)$
parameters and that  specializes to our stability condition for a triple when
$n=2$.
In fact he considers this notion for more general diagrams than the one above.
Presumably this stability condition governs, as for triples, the existence of
Hermitian
metrics on the bundles $E_i$ satisfying some generalized vortex equations
naturally
associated to the $(2n-1)$-tuple.

4. Our results have also been extended  in a different direction \cite{B-GP} to
parabolic
triples, that is to triples
in which the bundles are endowed with parabolic structures. The Higgs field can
be either a
 parabolic morphism  or a meromorphic morphism with simple poles at the
parabolic points and
whose residues respect the parabolic structure in some precise sense. In both
cases one can
prove a Hitchin-Kobayashi correspondence, although the metrics involved now
have
 singularities  at the parabolic points.

\vspace{12pt}

\noindent {\em Acknowledgements}.\ The authors would like to thank Alastair
King and Jun Li
for helpful conversations and  the following institutions
for their hospitality during the course of this project: The Mathematics
Institute of the
University of Warwick, England;   I.H.E.S., Fpance; the Mathelatics Department
of
UC Berkeley, USA; and C.I.M.A.T., Mexico.



\vspace{12pt}

\vspace{12pt}

\vspace{12pt}

\noindent Department of Mathematics, University of Illinois, 273 Altgeld Hall,
1409 W.
Green Street, Urbana, IL 61801, USA.

\vspace{12pt}

\noindent Universit\'e de Paris-Sud, Math\'ematiques, B\^atiment 425, 91405
Orsay, France.

\end{document}